\documentclass[10pt]{article}
\usepackage[T1]{fontenc} 
\usepackage[utf8]{inputenc}

\usepackage{lmodern}

\usepackage{enumerate}

\usepackage{lscape,graphicx,color,psfrag}
\usepackage{amssymb}
\usepackage[cmex10]{amsmath}
\usepackage{amsthm}
\usepackage{latexsym}
\usepackage{epsfig}
\usepackage{graphicx}

\DeclareMathOperator{\Ex}{\mathbb{E}}
\DeclareMathOperator{\Var}{Var}
\DeclareMathOperator{\sgn}{sgn}

\theoremstyle{plain}
\newtheorem{theorem}{Theorem}[section]
\newtheorem{corollary}[theorem]{Corollary}
\newtheorem{lemma}[theorem]{Lemma}
\newtheorem{proposition}[theorem]{Proposition}
\newtheorem{assumption}[theorem]{Assumption}

\theoremstyle{definition}

\theoremstyle{remark}
\newtheorem{remark}{Remark}[section]

\usepackage{chngcntr}
\usepackage{apptools}
\AtAppendix{\counterwithin{theorem}{section}}

\usepackage{appendix}

\usepackage[top=1in, bottom=1.25in, left=0.9in, right=0.9in]{geometry}

\usepackage[sort&compress,numbers]{natbib}
\usepackage[colorlinks=true,urlcolor=blue,citecolor=blue,linkcolor=blue]{hyperref}
\usepackage[noabbrev,nameinlink]{cleveref}     

\def\EMAIL#1{\href{mailto:#1}{#1}}

\begin{document}

\title{Solving Poisson's equation for birth--death chains: Structure, instability, and accurate approximation\footnote{Research supported by the Spanish Government under grants 
MICINN 
MTM2010-20808, MINECO ECO2015-66593-P, and   PID2019-109196GB-I00 / AEI / 10.13039/501100011033. The author has presented preliminary versions at the 7th International Conference on Network Games, Control and Optimization (NETGCOOP 2014, Trento, Italy) and the 29th European Conference on Operational Research (EURO 2018, Valencia, Spain).}}


\author{Jos\'e Ni\~no-Mora 
\\ Department of Statistics \\
    Carlos III University of Madrid \\
     28903 Getafe (Madrid), Spain \\  \EMAIL{jnimora@alum.mit.edu}, \href{http://alum.mit.edu/www/jnimora}{http://alum.mit.edu/www/jnimora} \\
      ORCID: \href{http://orcid.org/0000-0002-2172-3983}{0000-0002-2172-3983}}
 
\date{Published in \textit{Performance Evaluation}, vol.\ 145, 102163, 2021 \\ \vspace{.1in}
DOI: \href{https://doi.org/10.1016/j.peva.2020.102163}{10.1016/j.peva.2020.1021634}}

\maketitle

\begin{abstract}%
Poisson's equation plays a fundamental role as a tool for performance evaluation and optimization of Markov chains. For continuous-time birth--death chains with possibly unbounded transition and cost rates as addressed herein, when analytical solutions are unavailable its numerical solution can in theory be obtained by a simple forward recurrence. Yet, this may suffer from numerical instability, which can hide the structure of exact solutions. This paper presents three main contributions: (1) it establishes a structural result (convexity of the relative cost function) under mild conditions on transition and cost rates, which is relevant for proving structural properties of optimal policies in Markov decision models; (2) it elucidates the root cause, extent and prevalence of instability in numerical solutions by standard forward recurrence; and (3) it presents a novel forward--backward recurrence scheme to compute accurate numerical solutions. The results are applied to the accurate evaluation of the bias and the asymptotic variance, and are illustrated in an example.
\end{abstract}%

 
\newpage
\section{Introduction.}
\label{s:intro}
Consider a stochastic system whose state evolution is modeled by an ergodic continuous-time  birth--death Markov chain $\{X(t)\}_{t \geqslant 0}$ (see, e.g.,  \cite[Ch.\ 6]{ross19}) on the state space 
$\mathbb{N}_0 \triangleq  \{0, 1, 2, \ldots\}$, with birth and death rates 
 $\lambda_n > 0$ and $\mu_n > 0$ for $n \in \mathbb{N} \triangleq  \{1, 2, \ldots\}$, and $\lambda_0 > 0 = \mu_0$,  having steady-state probabilities $p_n$. 
Costs accrue at the state-dependent rates $c_n$, satisfying
\begin{equation}
\label{eq:zetak}
\sum_{n=0}^\infty |c_n| p_n < \infty,
\end{equation} 
so that the mean steady-state cost $\zeta \triangleq \sum_{n=0}^\infty c_n p_n$ is well defined and finite.

In applications to queueing systems, the state $n$ represents the number of jobs in the system, $\lambda_n$ is the state-dependent arrival rate, and $\mu_n$ can incorporate both service  and abandonment rates. The $c_n$ model a possibly nonlinear cost structure, which may account for different types of costs, such as holding or abandonment costs. See the examples in Remark \ref{re:pofzi}. There is an extensive literature on queueing models with costs, concerning their optimal dynamic control. See, e.g., \cite{stidh09} and references therein.

Birth, death and cost rates may all be \emph{unbounded}. Thus, the present setting is strictly more general than a discrete-time one, as analysis of a discrete-time Markov chain can be reduced to that of a continuous-time chain, whereas the reverse is only true for \emph{uniformizable} chains, having bounded transitions rates.

This paper addresses the exact and approximate numerical solutions to the \textit{Poisson equation}  for such a system, given by the second-order recurrence
\begin{equation}
\label{eq:pe0}
\lambda_0 b_{1} - \lambda_0 b_0  = \zeta-c_0, \quad 
\lambda_n b_{n+1} - (\lambda_n+ \mu_n) b_n + \mu_n b_{n-1}  = \zeta-c_n, \quad n \in \mathbb{N}.
\end{equation}
Note that (\ref{eq:pe0}) is a \textit{Poisson equation} in the usage of this term in the Markov chain literature, where it plays a fundamental role as a tool for performance evaluation and optimization. See, e.g., \cite{neveu71,nummelin91,glynnMeyn96}.
Ref.\ \cite[pp.\ 458--459]{meynTweed09}  surveys wide-ranging applications in applied probability, statistics and engineering, including  performance bounds,   analysis and   variance reduction in simulation, \emph{Markov decision processes} (MDPs) (see \cite{guoHL09}), and perturbation theory.
Note that the iterative solution of the dynamic programming optimality (Bellman) equations of an MDP model via Howard's \cite{howard60} \emph{policy improvement} method entails solving Poisson equations corresponding to stationary deterministic policies.

Concerning exact solutions, we aim to exploit explicit expressions to establish \emph{convexity} of functions $b\colon \mathbb{N}_0 \to \mathbb{R}$ solving (\ref{eq:pe0}), under mild conditions on transition and cost rates.
Note that a standard approach to prove structural properties of optimal policies in  MDP models relies on such a convexity property. See  \cite{blokSpieks17}. 

As for approximate numerical solutions, these are required when analytical solutions to (\ref{eq:pe0}) are unavailable. When $\zeta$ is not representable by a \emph{machine number}, due to finite-precision arithmetic, and is approximated by $\widehat{\zeta} \neq \zeta$, one may consider, ignoring other error sources, that the numerical recurrence actually solved is the 
modified Poisson equation (with $z = \widehat{\zeta}$)
\begin{equation}
\label{eq:pe1}
\lambda_0 b_{1} - \lambda_0 b_0  = z-c_0, \quad 
\lambda_n b_{n+1} - (\lambda_n+ \mu_n) b_n + \mu_n b_{n-1}  = z-c_n, \quad n \in \mathbb{N}.
\end{equation}

We will thus investigate  (\ref{eq:pe1}), viewing $z \in \mathbb{R}$ as an \emph{input parameter}, with the aim of elucidating how the approximation error in $\widehat{\zeta}$ is amplified into corresponding errors in the $\widehat{b}_n$ solving (\ref{eq:pe1}).

From the structure of (\ref{eq:pe1}), it is evident that there exists a solution $b$ for any choice of $z  \in \mathbb{R}$, which can be constructed by setting $b_0$ arbitrarily, and generating the remaining $b_n$ by the second-order linear forward recurrence
\begin{equation}
\label{eq:sonhlfrb}
 b_{1}  = \frac{z - c_0}{\lambda_0}  + b_0, \quad 
 b_{n+1} = 
\frac{z - c_n}{\lambda_n}  + \frac{\lambda_n+ \mu_n}{\lambda_n} b_n - \frac{\mu_n}{\lambda_n} b_{n-1}, \quad n = 1, 2, \ldots
\end{equation}
Clearly, for given $z$, the solutions $b$ to (\ref{eq:pe1}) are
unique up to an additive constant.

The present viewpoint contrasts with most work on Poisson's equation, which has focused on the choice $z = \zeta$ (see, e.g.,
\cite{glynnMeyn96}).  
A major reason for this is the following. Suppose that birth--death rates satisfy a \emph{Foster--Lyapunov drift condition} with \emph{weight function} $w\colon \mathbb{N}_0 \to [1, \infty)$, and let
$\mathbb{B}_w \triangleq \{f\colon \sup_{n} |f_n|/w_n < \infty\}$ be the Banach space  of  \emph{$w$-bounded} functions.
Then (see \cite[Theorem 7.1]{meynTweed93} and \cite[Prop.\ 7.11]{guoHL09}): 
(i) the chain satisfies a strong form of ergodicity (it is \emph{$w$-exponentially ergodic}); (ii) $\sum_{n=0}^\infty w_n p_n < \infty$; and (iii) for  any
$c \in \mathbb{B}_w$, (\ref{eq:pe1}), seen as an equation in $(z, b)$,  has a \emph{unique} solution in $\mathbb{R} \times \mathbb{B}_w$, up to an additive constant for $b$. Furthermore, this solution  is of the form $(z, b) = (\zeta, \beta+a)$ with $a \in \mathbb{R}$, where, writing as $\Ex_n[\cdot]$ the expectation starting from $X(0) = n$,
\begin{equation}
\label{eq:bias}
\beta_n \triangleq \Ex_n\bigg[\int_0^\infty (c_{ X(t)} - \zeta) \, dt\bigg]
\end{equation}
is the  \emph{bias} or \emph{relative cost} of starting from $n$ rather than from steady state.
Note that $\beta$ is characterized as the only $b \in \mathbb{B}_w$ solving  (\ref{eq:pe0}) for which  
 \begin{equation}
\label{eq:adc1}
\sum_{n =0}^\infty b_n p_n = 0.
\end{equation}

As for interpretation, if $b$ with $b_m = 0$ solves (\ref{eq:pe0}),  then $b_n$ is the \emph{relative cost} of starting from $n$ rather than $m$, so    $b_n = \Ex[\int_0^\infty (c_{X^n(t)} - c_{X^m(t)}) \, dt]$, for realizations  $\{X^n(t)\}_{t \geqslant 0}$ and $\{X^m(t)\}_{t \geqslant 0}$ of the chain  starting from $n$ and $m$.

We further consider the function $\varphi\colon \mathbb{N}_0 \to \mathbb{R}$ defined by
\begin{equation}
\label{eq:varphidef}
\varphi_n \triangleq \beta_{n+1} - \beta_n = b_{n+1} - b_n, \quad n \in \mathbb{N}_0,
\end{equation}
for any $b$ solving (\ref{eq:pe0}).
We will refer to $\varphi$ as the \emph{marginal relative cost} function, as $\varphi_n$ is the relative cost of starting from $n+1$ rather than $n$.
Note that $\varphi$ is determined by the following first-order linear recurrence, which reformulates (\ref{eq:pe0}):
\begin{equation}
\label{eq:varphirec}
\lambda \varphi_0  = \zeta,  \quad  \lambda \varphi_n - \mu_n 
\varphi_{n-1} = \zeta - c_n, \quad  n \in \mathbb{N}.
\end{equation}
Similarly to (\ref{eq:pe1}), when (\ref{eq:varphirec}) is solved numerically substituting $z = \widehat{\zeta}$ for $\zeta$, the recurrence that is actually solved (ignoring other sources of error) is
\begin{equation}
\label{eq:pe2}
\lambda f_0  = z,  \quad  \lambda f_n - \mu_n 
f_{n-1} = z - c_n, \quad  n \in \mathbb{N}.
\end{equation}

To illustrate, consider (cf.\ \cite[\S 4.2]{nmcor19}), the M/M/$1$+M  queue with deadlines to the end of service, with $\lambda_n \triangleq \lambda$, 
$\mu_n \triangleq  \mu + n \theta$ and $c_n \triangleq n \theta$, where $\lambda, \mu, \theta > 0$ are the arrival, service and abandonment rates.
From (\ref{eq:pe2}), we obtain the recursion
\begin{equation}
\label{eq:pepiim}
f_0  = \frac{z}{\lambda},  \quad  f_n  = \frac{z - n \theta}{\lambda} + \frac{\mu + n \theta}{\lambda}
f_{n-1}, \quad  n = 1, 2, \ldots
\end{equation}

Consider now the instance with $\lambda =  0.9$, $\mu = 1$ and $\theta = 0.5$, for which 
\begin{equation}
\label{eq:zetaex}
\zeta = \frac{1}{10} \bigg(\frac{81}{5 e^{9/5}-14}-1\bigg).
\end{equation}
 Computing (\ref{eq:zetaex}) with MATLAB gives the double-precision floating-point number (see  \cite[Ch.\ 2]{higham02}) $\widehat{\zeta} \approx 0.398515613690624$. 
Table \ref{tab:numex1} shows, in the $\widehat{\varphi}_{n}$ column, the first thirty  $f_n$ computed through (\ref{eq:pepiim}) using $z = \widehat{\zeta}$.
After growing to $1$ at first, $\widehat{\varphi}_{n}$  diverges to \emph{minus} infinity.
The same behavior results when $z$ is set to the next larger machine number, $\widehat{\zeta} + 2^{-54}$.
Yet, when the following machine number, $z = \widetilde{\zeta} = \widehat{\zeta} + 2^{-53}$, is used, the resulting $f_n$, written as $\widetilde{\varphi}_{n}$,  eventually diverge to \emph{plus} infinity.
Thus, approximate computation of the $\varphi_n$ through (\ref{eq:pepiim}), and hence of the $b_n$, suffers from an \emph{explosive numerical instability} with respect to unavoidable errors in the approximation to the input $\zeta$.

Experimentation with other models reveals that such instability is not exceptional.
Thus, e.g., \cite{steckHend07} reports numerical instabilities preventing accurate numerical solution of Poisson's equation in the queueing models considered there.

Table \ref{tab:numex1} further shows the probabilities $p_n$, to gain insight on the relation between these and the magnitudes of errors in the computed approximations to $\varphi_n$. The loss of significant digits of accuracy in the latter is evident only for unlikely states, with  such inaccuracies growing steeply  as the $p_n$ get  smaller. 
Note that the accurate estimation of performance metrics associated to very low-probability states, corresponding to \emph{rare events}, is   of considerable interest in a variety of areas, including computer-communication systems, reliability, finance, etc.,  and is the subject of major research attention. 
See, e.g.,  the survey \cite{junejaShah06}.

\begin{table}[tbh!]
\caption{Approximate numerical computation of $\varphi_n = b_{n+1}-b_n$: explosive instability.} \label{tab:numex1}
\scriptsize 
 \vspace{.1in}
\begin{center}
\begin{tabular}{|rlll||llll|}  
\multicolumn{1}{c}{$n$} & \multicolumn{1}{c}{$p_n$} & \multicolumn{1}{c}{$\widehat{\varphi}_{n}$} & \multicolumn{1}{c}{$\widetilde{\varphi}_{n}$} & \multicolumn{1}{c}{$n$} & \multicolumn{1}{c}{$p_n$} & \multicolumn{1}{c}{$\widehat{\varphi}_{n}$} & \multicolumn{1}{c}{$\widetilde{\varphi}_{n}$}\\ \hline
$0$ & $5.0 \times 10^{-1}$ & $0.44279513$ & $0.44279513$ &   $15$ & $1.9 \times 10^{-11}$ & $0.9388453$ & $0.9388507$\\
$1$ & $3.0\times 10^{-1}$ & $0.62523145$ & $0.62523145$ & $16$ & $1.9 \times 10^{-12}$ & $0.9423595$ & $0.9424134$ \\
$2$ & $1.4 \times 10^{-1}$ & $0.72108723$ & $0.72108723$ & $17$ & $1.8 \times 10^{-13}$ & $0.9454791$ & $0.9460477$ \\
$3$ & $4.9 \times 10^{-2}$ & $0.77914855$ & $0.77914855$ & $18$ & $1.6 \times 10^{-14}$ & $0.9481179$ & $0.9544357$ \\
$4$ & $1.5\times 10^{-2}$ & $0.81773474$ & $0.81773474$ & $19$ & $1.4 \times 10^{-15}$ & $0.9486149$ & $1.0223229$\\
$5$ & $3.7 \times 10^{-3}$ & $0.84509690$ & $0.84509690$ &  $20$ & $1.1 \times 10^{-16}$ & $0.9258667$ & $1.8267420$ \\
$6$ & $8.4 \times 10^{-4}$ & $0.86544800$ & $0.86544800$ &  $21$ & $8.9 \times 10^{-18}$ & $0.6066468$ & $1.2 \times 10^1$\\
$7$ & $1.7 \times 10^{-4}$ & $0.88114626$ & $0.88114626$ &  $22$ & $6.6 \times 10^{-19}$ & $-3.7 \times 10^0$ & $1.5 \times 10^2$\\
$8$ & $3.0 \times 10^{-5}$ & $0.89360768$ & $0.89360768$ &  $23$ & $4.8 \times 10^{-20}$ & $-6.4 \times 10^1$ & $2.1 \times 10^3$\\
$9$ & $5.0 \times 10^{-6}$ & $0.90373094$ & $0.90373094$ &  $24$ & $3.3 \times 10^{-21}$ & $-9.3 \times 10^2$ & $3.0 \times 10^4$\\
$10$ & $7.4 \times 10^{-7}$ & $0.91211248$ & $0.91211248$ & $25$ & $2.2 \times 10^{-22}$ & $-1.4 \times 10^4$ & $4.5 \times 10^5$ \\
$11$ & $1.0 \times 10^{-7}$ & $0.91916304$ & $0.91916304$ &  $26$ & $1.4 \times 10^{-23}$ & $-2.2 \times 10^5$ & $7.0 \times 10^6$\\
$12$ & $1.3 \times 10^{-8}$ & $0.92517434$ & $0.92517435$ &  $27$ & $8.8 \times 10^{-25}$ & $-3.5 \times 10^6$ & $1.1 \times 10^8$\\
$13$ & $1.6 \times 10^{-9}$ & $0.93035909$ & $0.93035915$ &  $28$ & $5.3 \times 10^{-26}$ & $-5.8 \times 10^7$ & $1.9 \times 10^9$\\
$14$ & $1.8 \times 10^{-10}$ & $0.93487590$ & $0.93487647$ &  $29$ & $3.1 \times 10^{-27}$ & $-1.0 \times 10^9$ & $3.2 \times 10^{10}$\\ \hline
\end{tabular}
\end{center}
\end{table}

Regarding marginal relative costs $\varphi_n$, they are of direct interest in a major application of Poisson's equation: the \emph{one-step policy improvement} (OSPI) method to the design of  scalable heuristic policies for certain multidimensional MDPs, e.g., those concerning the optimal routing  of a job stream to parallel service stations,  each with its own queue. 
One starts with a tractable \emph{static policy} (state-independent)  under which the queues evolve independently.
Then, a single step of  Howard's \emph{policy improvement} algorithm (see \cite{howard60}) for MDPs  gives a \emph{dynamic policy} (state-dependent) with better cost performance.

The latter is an \emph{index policy}, where the \emph{index}  for each station is precisely the marginal relative cost $\varphi_n$, and each arrival is  routed to a station of currently lowest index value. 
See, e.g., \cite{krish90,bhulai06,hyytiaRight16,nmcor12,nmejor12,nmcor19,tijms17,bhulai17}.

Intuition would suggest that, if  such a policy routes an arrival to a station in a given system state, it should prescribe the same action if that station was less loaded, other things being equal.
This would be ensured if one could prove that
$\varphi_n$  is  \emph{nondecreasing} in $n$ or, equivalently, that $\beta_n$  is \emph{convex}.  
Note that, in Table \ref{tab:numex1}, the computed $\varphi_n$ increase at first, but, when numerical instabilities set in, such a monotonicity is lost, raising doubts on the behavior of the exact $\varphi_n$. 

A related metric of interest is the \emph{asymptotic variance} as $t \to \infty$ of the average cost up to time $t$, $\bar{C}(t) \triangleq (1/t) \int_0^t c_{ X(s)} \, ds$, which is defined by
\[
\sigma^2 \triangleq \lim_{t \to \infty} \, t \, \Var\big[\bar{C}(t)\big].
\]
Under the aforementioned  drift condition, and provided that $c^2 \in \mathbb{B}_w$,   
$\sigma^2$ is well defined and finite, being given by (see \cite[Theorem 4.4]{glynnMeyn96})
\begin{equation}
\label{eq:sigma22}
\sigma^2 = 2 \sum_{n =0}^\infty b_n (c_n - \zeta) p_n = 2 \sum_{n =0}^\infty \beta_n c_n   p_n
\end{equation}
for any $b \in \mathbb{B}_w$ solving (\ref{eq:pe0}), 
and $\bar{C}(t)$ satisfies the \emph{functional central limit theorem} 
\[
\sqrt{t} \, \frac{\bar{C}(t) - \zeta}{\sigma} \Longrightarrow \mathrm{N}(0, 1) \textup{  as  } t \to \infty,
\]
where $\Longrightarrow$ denotes weak convergence and $\mathrm{N}(0, 1)$ is the standard normal distribution.
This result can be used, e.g.,  to set simulation run lengths to obtain  confidence intervals for $\zeta$.
See, e.g.,   \cite{whittAFMP92,steckHend07,niHend15}.

\subsection{Goals and contributions}
\label{s:gac}
We address the following research goals: 1)
identify mild conditions on transition and cost rates ensuring that the $\varphi$ solving (\ref{eq:pe1}) is  nondecreasing (equivalently, the $b$'s solving (\ref{eq:pe0}) are convex);
2) elucidate the  root cause, extent and prevalence of numerical instability in computed solutions to Poisson's equation due to substituting approximations $\widehat{\zeta}$ for $\zeta$;  and 
3) obtain means of  computing accurate approximate solutions when analytical solutions are not available.

Regarding the first goal, consider the following assumption, where 
we write  
$d_n \triangleq \mu_n - \lambda_n$ and 
use the backward difference notation $\Delta x_n \triangleq x_n - x_{n-1}$.

\begin{assumption}
\label{ass:hmulambda} 
\begin{itemize}
\item[\textup{(i)}] $d$ is 
\begin{itemize}
\item[\textup{(i.a)}]
nondecreasing $(\Delta d_{n} \geqslant 0),$ with $\Delta d_1 > 0;$ and
\item[\textup{(i.b)}] concave $(\Delta d_{n+1} \leqslant \Delta d_{n}).$ 
\end{itemize}
\item[\textup{(ii)}] $c$ is 
\begin{itemize}
\item[\textup{(ii.a)}] nonnegative $(c_n \geqslant 0)$ and nondecreasing $(\Delta c_n \geqslant 0);$ and
\item[\textup{(ii.b)}] convex $(\Delta c_{n+1} \geqslant \Delta c_{n}).$
\end{itemize}
\end{itemize}
\end{assumption}

We have the following result. 

\begin{theorem}
\label{the:pofzi} Under Assumption $\ref{ass:hmulambda},$
$\varphi$ is nondecreasing $(\beta$ is convex$).$ 
\end{theorem}

\begin{remark}
\label{re:pofzi} 
\begin{itemize}
\item[\textup{(a)}] Assumption \ref{ass:hmulambda} is mild and widely satisfied in birth--death queueing models. 
Thus, consider the following broad $m$-server model with possible customer balking and abandonment, which encompasses a variety of standard models:
$\lambda_n = \lambda (1-\alpha_n)$, with $0 \leqslant \alpha_n < 1$ the state-dependent balking probability, $\mu_n = \min(m, n) \mu + g(m, n) \theta$, where $\mu > 0$ and $\theta \geqslant 0$ are the service and abandonment rates, and $g(m, n) \triangleq  (n-m)^+$ if customers can only abandon prior to entering service, while $g(m, n) \triangleq  n$ otherwise. Further, $c_n = c^{\textup{ab}} g(m, n) \theta + c_n^{\textup{h}}$, where $c^{\textup{ab}} \geqslant 0$ and $c_n^{\textup{h}} \geqslant 0$ are the cost per abandonment and the holding cost rate. It is readily verified that this model satisfies Assumption \ref{ass:hmulambda} if the following holds: (i) $\alpha_n$ is nondecreasing and concave; (ii) in the case that customers can only abandon before entering service, $\theta \leqslant \mu$, i.e., the mean time to abandon ($1/\theta$) is not shorter than the mean service time ($1/\mu$); and (iii)  $c_n^{\textup{h}}$ is nondecreasing and convex. These are all intuitively sound conditions.
\item[\textup{(b)}] Theorem \ref{the:pofzi} is also relevant for birth--death models arising in population biology. Thus, consider the classic linear birth--death model with immigration, with $\lambda_n = n \lambda + \alpha$ and $\mu_n = n \mu$. In the case of a population whose size is costly, e.g., a pest population as in \cite{kyriDim16}, one may take $c_n$ to be the cost of having $n$ individuals. Then, Assumption \ref{ass:hmulambda} clearly holds provided that $c$ is nonnegative, nondecreasing and convex.
\item[\textup{(c)}] Note that the conditions in Assumption \ref{ass:hmulambda} were first formulated in \cite[Ass.\ 7.1]{nmmp02}, though with a different purpose: they were shown there to be sufficient conditions for existence of the \emph{Whittle index} characterizing optimal policies in a broad birth--death admission control model.
\end{itemize}
\end{remark}

As for the second goal, we explicitly identify and analyze the \emph{error amplification factors} (see, e.g., \cite{gautschi61} for early use of such a concept) that characterize how the error in the approximation $\widehat{\zeta}$ to the input $\zeta$ propagates, in the standard \emph{forward recurrence} scheme,  to produce errors in the computed approximations $\widehat{b}_n$ and $\widehat{\varphi}_n$ to the
$b_n$ and $\varphi_n$ solving (\ref{eq:pe0}) and (\ref{eq:varphirec}).
See Proposition \ref{pro:phihatsol}.
Such results are further used to analyze the accuracy of computed approximations to the bias $\beta$ and the asymptotic variance $\sigma^2$. See Proposition \ref{pro:nbetasigma2phi}.

Concerning the third goal, we analyze the approximation errors resulting from a \emph{backward recurrence} scheme, which has not been previously considered for this model and has substantially improved accuracy for large states. See Propositions \ref{pro:brphihatsol} and \ref{pro:compfrbr}. 
We further propose a novel mixed \emph{forward--backward recurrence} scheme that outperforms both forward and backward recurrence  in terms of accuracy. 
See Propositions \ref{pro:fbrphihatsol} and \ref{pro:fbnbetasigma2phi}.
The effectiveness of the proposed approach is demonstrated for the motivating example above. See \S\ref{s:anex}.

\subsection{Organization of the paper}
\label{s:ootp}
The remainder of the paper proceeds as follows.
\S \ref{s:rw} reviews related work.
\S \ref{s:bdma} presents results  on birth--death chains, some of them new, that are required for subsequent analyses.
\S  \ref{s:ersp} gives expressions for the exact solution to Poisson's equation, as well as new expressions for the bias and the asymptotic variance in terms of marginal relative costs; it further gives the proof of Theorem \ref{the:pofzi}.
\S \ref{s:fraeas} develops an error  analysis of computed solutions to Poisson's equation through forward recurrence.
\S \ref{s:abr} analyzes a backward recurrence scheme with improved accuracy for large states.
\S \ref{s:afbr} presents a mixed forward--backward recurrence scheme, with improved accuracy with respect to the pure schemes.
\S \ref{s:anex} illustrates the results in an example.
\S \ref{s:cr} concludes. Appendix
\ref{a:pofzi} lays the groundwork for the proof of Theorem \ref{the:pofzi}.

Note that an early version of this work appeared in the proceedings \cite{nmnetgcoop14}.

\section{Further related work}
\label{s:rw}
In addition to the work referred to in \S\ref{s:intro}, we  briefly review next two further related strands of research.
The first one focuses on the Poisson equation for discrete-time Markov chains. \cite{makoShw02} reviews structural results from a probabilistic interpretation viewpoint, while 
\cite{bhulaiSpieks03} establishes uniqueness of solutions to Poisson's equation in such a setting. 
Another line of work aims to elucidate the structure of solutions to Poisson's equation for discrete-time birth--death chains and extensions. See, e.g.  \cite{dendetal13,binietal16}.

Another closely related relevant stream of work is the 
wide literature on numerical stability of computed solutions to linear recurrences, both in general and applied to the evaluation of special functions. In a classic study, Gautschi addressed in 
\cite{gautschi61}  the recursive numerical solution of a  general first-order linear recurrence relation, assuming that the only source of error is the approximation of the initial condition.
He considered the propagation of such an error to the computation of successive terms, as measured by relative \emph{error amplification factors}, under two recursive schemes, standard \emph{forward recurrence} and \emph{backward recurrence}. The paper, which refers to earlier work on the effective use of the latter type of recurrence, elucidates in which cases each scheme is more accurate.
Note however that, although (\ref{eq:varphirec}) is a first-order linear recurrence, the analysis in \cite{gautschi61} does not directly apply to it, as that paper considers that the source of error is in the initial condition, whereas here we consider that it lies in the approximation to $\zeta$ in the right-hand side of the equation.  

Early work on second-order linear recurrence relations, regarding  
``\emph{the possibility and the prevention of numerical instability}'', is reviewed in \cite{gautschi67}, which focuses on the homogeneous case, while, e.g., 
\cite{olver67,amosBurg73} consider the non-homogeneous case. In such papers, numerically stable schemes, including backward recurrence, for accurately computing so-called \emph{minimal solutions} are developed.
Yet, again, such work does not directly apply to the analysis of Poisson's equation (\ref{eq:pe1}), since the source of error addressed here is in the right-hand sides, which is not considered in the aforementioned papers.

For more recent work see, e.g., \cite{barrioetal03}, which develops 
rounding-error bounds for the numerical solution of higher-order linear recurrence relations. 

Even though, for the reasons outlined above, such work does not appear to be directly applicable to the recurrence relations herein, still, key ideas developed in that field such as 
error amplification factor analysis, and the use of backward recurrence to overcome numerical instability, play a central role in this paper.

Yet, such general ideas need to be adapted to the present setting, by exploiting the rich special structure and the probabilistic interpretation of the recurrences herein.
To the best of the author's knowledge, such an analysis has not been undertaken before, and is  hence a novel contribution of this work. 

\section{Required results on birth--death Markov chains}
\label{s:bdma}
This section presents results that are required for the ensuing analyses, mostly known, but also some new extensions (which the author has not found in the literature), on mean first-passage times and costs for birth--death chains. 

Consider a birth--death process $\{X(t)\}_{t \geqslant 0}$ with costs as outlined in \S \ref{s:intro}.
Its \emph{potential coefficients} $\pi_n$ are the unnormalized state probabilities, 
\[
\pi_0 \triangleq 1, \quad
\pi_n \triangleq \frac{\lambda_0 \cdots \lambda_{n-1}}{\mu_1 \cdots \mu_{n}}, \quad n \in \mathbb{N},
\] 
which satisfy the recurrence
\begin{equation}
\label{eq:rhoki}
\pi_{0} = 1, \quad 
\lambda_{n-1} \pi_{n-1} = \mu_n \pi_{n}, \quad n \in \mathbb{N}.
\end{equation}

We assume that 
\begin{equation}
\label{eq:uniqdbr}
\sum_{n=0}^\infty \frac{1}{\lambda_n \pi_n} = \infty \quad
\textup{and} \quad
\sum_{n=0}^\infty \pi_n < \infty,
\end{equation}
which (see \cite[Theorem 2(a)]{karlMcG57}) is a necessary and sufficient condition for 
the process to be  uniquely determined by its birth--death rates and ergodic.
The following is
a standard sufficient condition for (\ref{eq:uniqdbr}), where $\rho_n \triangleq \lambda_n/\mu_n$:
\begin{equation}
\label{eq:lmunrho}
\limsup_{n \to \infty} \, \rho_n < 1.
\end{equation}

The steady-state probabilities are given by
\begin{equation}
\label{eq:pirhoi}
p_n = \frac{\pi_n}{\sum_{j=0}^\infty \pi_{j}} = \pi_n p_{0}, \quad n \in \mathbb{N}_0.
\end{equation}

We will consider the functions 
\begin{equation}
\label{eq:PC}
P_n \triangleq \sum_{j=0}^n p_{j}, \quad
C_n \triangleq \sum_{j=0}^n c_{j} p_{j}  \quad \textup{and} \quad  Z_n \triangleq \frac{C_n}{P_n}, \quad n \in \mathbb{N}_0,
\end{equation}
which satisfy
\begin{equation}
\label{eq:PCZ1is}
P_n \to 1, \quad C_n \to \zeta  \quad \textup{and} \quad  Z_n \to \zeta \textup{ as} \quad  n \to \infty.
\end{equation}

Thus, $P_n$ is the cumulative distribution function of the $p_n$. We also consider the corresponding \emph{tail} functions defined by
\begin{equation}
\label{eq:PCtails}
\bar{P}_n \triangleq 1-P_n, \quad
\bar{C}_n \triangleq \zeta-C_n  \quad \textup{and} \quad  \bar{Z}_n \triangleq \zeta-Z_n, \quad n \in \mathbb{N}_0.
\end{equation}

Letting $\tau_n \triangleq \min \{t \geqslant 0\colon X(t) = n\}$ be the \emph{first-passage time to state $n$}, let 
$T_n^+ \triangleq \Ex_n[\tau_{n+1}]$ and $T_n^- \triangleq \Ex_n[\tau_{n-1}]$ be the mean first-passage times from $n$ to $n+1$, and from $n$ to $n-1$, and write as $H_n^+ \triangleq \Ex_n[\int_0^{\tau_{n+1}} c_{X(s)} \, ds]$ and 
$H_n^- \triangleq \Ex_n[\int_0^{\tau_{n-1}} c_{X(s)} \, ds]$ the corresponding mean costs.
We further define
\begin{equation}
\label{eq:Znpm}
Z_n^+ \triangleq \frac{H_n^+}{T_n^+} \quad \textup{and} \quad
Z_{n+1}^- \triangleq \frac{H_{n+1}^-}{T_{n+1}^-}, \quad n \in \mathbb{N}_0.
\end{equation}

The following recurrences follow from elementary probabilistic arguments:

\begin{equation}
\label{eq:Pi}
\lambda_0 T_0^+ = 1, \quad 
\lambda_n T_n^+ - \mu_n T_{n-1}^+ = 1, \quad n \in \mathbb{N},
\end{equation}
which is given in \cite[Eq.\ (1.5a)]{keilson65} (see its derivation in \cite[Ch.\ 6.3]{ross19})
and 
\begin{equation}
\label{eq:Hi}
\lambda_0 H_0^+ = c_0, \quad \lambda_n H_n^+ - \mu_n H_{n-1}^+ = c_n, \quad n \in \mathbb{N}.
\end{equation}

In \cite[Eq.\ (1.8)]{keilson65}, it is shown  that $T_n^+$ and $P_n$ are linked by
\begin{equation}
\label{eq:Tn}
T_n^+ =  \frac{P_n}{\lambda_n p_n},
\end{equation}
and the same argument given there yields that
\begin{equation}
\label{eq:Hn}
H_n^+ =  \frac{C_n}{\lambda_n p_n}.
\end{equation}

Note that it follows immediately from (\ref{eq:Tn}) that, as $n \to \infty,$
\begin{equation}
\label{eq:Tncons}
T_n^+ \to \infty \enspace \textup{if and only if} \enspace \lambda_n p_n \to 0,
\end{equation}
and 
\begin{equation}
\label{eq:supTncons}
\sup_n \, T_n^+ = \infty \enspace \textup{if and only if} \enspace \inf_n \, \lambda_n p_n = 0,
\end{equation}
so the $T_n^+$ are bounded if and only if the $\lambda_n p_n$ are bounded away from $0$.

We next give a sufficient condition for $T_n^+  \to \infty$
in terms of the $\rho_n$ in (\ref{eq:lmunrho}).

\begin{lemma}
\label{lma:tnpinfty}
Under condition $(\ref{eq:lmunrho}),$ $T_n^+ \to \infty$ as $n \to \infty.$
\end{lemma}
\begin{proof}
Write $\rho^* \triangleq \limsup_{n \to \infty} \, \rho_n$. 
Pick $\varepsilon > 0$  such that $\rho^* + \varepsilon < 1$, and let $N$ be such that $\rho_n < \rho^* + \varepsilon$ for $n > N$. Then 
\[
T_n^+ = \frac{1}{\lambda_n} + \frac{T_{n-1}^+}{\rho_n} > \frac{T_{n-1}^+}{\rho^* + \varepsilon}, \quad 
n > N,
\]
from which it follows that $T_n^+ \to \infty$ as $n \to \infty.$
\end{proof}

The following recurrences are also easily obtained:
\begin{equation}
\label{eq:Pim}
\mu_n T_n^- - \lambda_n T_{n+1}^- = 1, \quad n \in \mathbb{N},
\end{equation}
which is given in \cite[Eq.\  (1.10)]{keilson65}, and its counterpart for costs, 
\begin{equation}
\label{eq:Him}
\mu_n H_n^-  - \lambda_n H_{n+1}^- = c_n, \quad n \in \mathbb{N}.
\end{equation}
Yet, unlike recurrences (\ref{eq:Pi}) and (\ref{eq:Hi}), which determine the $T_n^+$ and $H_n^+$, (\ref{eq:Pim}) and (\ref{eq:Him}) do not determine the $T_n^-$ and $H_n^-$, as they lack boundary conditions.

These are provided by the case $n = 0$ in the following result, which gives 
analogous relations to (\ref{eq:Tn}) and (\ref{eq:Hn}) linking $T_{n+1}^-$ to $\bar{P}_n$ and $H_{n+1}^-$ to $\bar{C}_n$.
Note that the identity in Lemma \ref{lma:CPm}(a) is given in \cite[Eq.\ (1.11)]{keilson65}. Yet, the argument outlined there bypasses a critical detail, as it entails summation of an infinite telescoping series, which requires showing that $\lambda_n p_n T_{n+1}^- \to 0$ as $n \to \infty$. 
We avoid this in the proof below by using instead an induction argument, drawing on the theory of renewal reward processes (see, e.g., \cite[Ch.\  7]{ross19}).

\begin{lemma}
\label{lma:CPm} For $n \in \mathbb{N}_0,$
\begin{itemize}
\item[\textup{(a)}]
$\displaystyle 
T_{n+1}^- = \frac{\bar{P}_n}{\lambda_n p_n};$
\item[\textup{(b)}] $\displaystyle 
H_{n+1}^- = \frac{\bar{C}_n}{\lambda_n p_n}.
$ 
\end{itemize}
\end{lemma}
\begin{proof}
(a) 
We use induction. For $n = 0$, consider the regenerative cycle starting from state  $0$ until the first return to $0$, and denote by $T_{00}$ its mean duration. By the \emph{renewal reward theorem}, and using that $T_{00} = 1/\lambda_0 + T_1^-$, we have
\[
p_0 = \frac{\Ex[\textup{cost during a cycle}]}{\Ex[\textup{duration of a cycle}]} = \frac{1/\lambda_0}{1/\lambda_0 + T_1^-},
\]
and hence $T_{1}^- = \bar{P}_0/\lambda_0 p_0$.
Suppose now that the result holds for $n-1$, so that $T_{n}^- = \bar{P}_{n-1}/\lambda_{n-1} p_{n-1}$. Then, using this and  (\ref{eq:Pim}) we obtain
\[
\lambda_n p_n T_{n+1}^- = \mu_n p_n T_n^- - p_n = 
\lambda_{n-1} p_{n-1} T_n^- - p_n = \bar{P}_{n-1} - p_n = \bar{P}_{n},
\]
which completes the induction.

(b) For $n = 0$,  let $H_{00}$ be the mean cost accrued over the cycle in part (a). By the renewal reward theorem, and since $H_{00} = c_0/\lambda_0 + H_1^-$, we have
\[
\zeta = \frac{\Ex[\textup{cost during a cycle}]}{\Ex[\textup{duration of a cycle}]} = \frac{c_0/\lambda_0 + H_1^-}{1/\lambda_0 + T_1^-},
\]
whence $H_{1}^- = \bar{C}_0/\lambda_0 p_0$.
Suppose now the result holds for $n-1$, so that $H_{n}^- = \bar{C}_{n-1}/\lambda_{n-1} p_{n-1}$. Then, using this and  (\ref{eq:Him}) gives, as required.
\[
\lambda_n p_n H_{n+1}^- = \mu_n p_n H_n^- - c_n p_n = 
\lambda_{n-1} p_{n-1} H_n^- - c_n p_n = \bar{C}_{n-1} - c_n p_n = \bar{C}_{n}.
\]
\end{proof}

\begin{remark}
\label{re:CP}
\begin{itemize}
\item[\textup{(a)}] From (\ref{eq:Tn}), (\ref{eq:Hn}) and Lemma \ref{lma:CPm} we obtain
 \begin{equation}
\label{eq:Tnmkeil}
P_n T_{n+1}^- = 
\bar{P}_n T_n^+ \quad \textup{and} \quad 
C_n H_{n+1}^- = 
\bar{C}_n H_n^+.
\end{equation} 
\item[\textup{(b)}] It follows from  (\ref{eq:Tn}) and (\ref{eq:Hn})  that the $Z_n$ in (\ref{eq:PC}) and $Z_n^+$ in (\ref{eq:Znpm}) satisfy
\begin{equation}
\label{eq:Zirel}
Z_n^+ = Z_n. 
\end{equation}
\item[\textup{(c)}] As for the $Z_{n+1}^-$ in (\ref{eq:Znpm}), 
from  Lemma \ref{lma:CPm} and (\ref{eq:ei2}) we can write (see (\ref{eq:PCtails}))
\begin{equation}
\label{eq:Zirelm}
Z_{n+1}^- = \frac{\bar{C}_n}{\bar{P}_n} = 
\frac{\zeta - C_n}{1-P_n} = Z_n + \frac{\bar{Z}_n}{\bar{P}_n} = \zeta + P_n \frac{\bar{Z}_n}{\bar{P}_n}.
\end{equation}
Thus, $Z_{n+1}^- - \zeta$ is asymptotically equivalent to the tail ratio ${\bar{Z}_n}/{\bar{P}_n}$.  Note that, typically, $Z_{n+1}^-$ will \emph{not} converge to $\zeta$ as $n \to \infty$. 
\end{itemize}
\end{remark}

Now, the aforementioned result that $\lambda_n p_n T_{n+1}^- \to 0$ as $n \to \infty$, and the corresponding result for $H_{n+1}^-$, follow immediately from Lemma \ref{lma:CPm} and (\ref{eq:PCZ1is}).
 
\begin{corollary}
\label{cor:lnpntnp1m}
As $n \to \infty,$ 
\begin{itemize}
\item[\textup{(a)}]
$\displaystyle 
\lambda_n p_n T_{n+1}^- \to 0;$
\item[\textup{(b)}] $\displaystyle 
\lambda_n p_n H_{n+1}^- \to 0.
$ 
\end{itemize}
\end{corollary}

We will also consider $T_{0n} \triangleq \Ex_0[\tau_n]$, the \emph{mean first-passage time from $0$ to $n$},  
$T_{n0} \triangleq \Ex_n[\tau_0]$, the \emph{mean first-passage time from $n$ to $0$}, and the corresponding mean costs $H_{0n} \triangleq \Ex_0[\int_0^{\tau_n} c_{X(s)} \, ds]$ and $H_{n0} \triangleq \Ex_n[\int_0^{\tau_0} c_{X(s)} \, ds].$
 Note that
\begin{equation}
\label{eq:T0np1H0np1}
 T_{0n} = \sum_{j=0}^{n-1} T_j^+ = \sum_{j=0}^{n-1} \frac{P_j}{\lambda_{j} p_{j}}
 \quad
 \textup{and}
 \quad 
 H_{0n} = \sum_{j=0}^{n-1} H_j^+ = \sum_{j=0}^{n-1} \frac{C_j}{\lambda_{j} p_{j}}, 
 \end{equation}
 whereas
 \begin{equation}
\label{eq:Tn0}
T_{n0} = \sum_{j=0}^{n-1} T_{j+1}^- = \sum_{j=0}^{n-1} \frac{\bar{P}_j}{\lambda_{j} p_{j}}
\quad
 \textup{and}
 \quad 
 H_{n0} = \sum_{j=0}^{n-1} H_{j+1}^- = \sum_{j=0}^{n-1} \frac{\bar{C}_j}{\lambda_{j} p_{j}}.
\end{equation} 
 
\begin{remark}
\label{re:T0np1div}
\begin{enumerate}[(a)]
\item
From (\ref{eq:Tn}) and (\ref{eq:T0np1H0np1}) we obtain
\[
\lim_{n \to \infty}\, T_{0n} = \sum_{j=0}^\infty T_j^+ = \sum_{j=0}^\infty \frac{P_j}{\lambda_j p_j} =  \infty,
\]
since the last series diverges by the limit comparison test (since $P_j \to 1$ as $j \to \infty$) using that $\sum_{j=0}^\infty 1/\lambda_j p_j = \infty$ by 
(\ref{eq:uniqdbr}).
\item We will refer to $T_{\infty 0} \triangleq \lim_{n \to \infty} T_{n0} = \sum_{j=1}^\infty T^-_j$, which is given by
\begin{equation}
\label{eq:tinfty0}
T_{\infty 0} 
= \sum_{n=0}^\infty \frac{\bar{P}_{n}}{\lambda_{n} p_{n}}.
\end{equation}
In Feller's \cite{feller59} classic boundary classification for birth--death processes, 
if $T_{\infty 0} < \infty$ the process is said to have an \emph{entrance boundary at  infinity}, and is then  exponentially ergodic. See \cite[Theorem 8.1]{CallKeil73II}. Otherwise, the process is said to have a \emph{natural boundary at infinity}. From \cite[Eq.\ (6.5) and Theorem 6.4(ii)]{CallKeil73I} it follows that, under condition (\ref{eq:lmunrho}),
\begin{equation}
\label{eq:tinf0fcond}
T_{\infty 0} < \infty \enspace \textup{if and only if} \enspace \sum_{n=1}^\infty \frac{1}{\mu_n} < \infty.
\end{equation} 
Note that, in most applied models, the rightmost series in  (\ref{eq:tinf0fcond}) diverges. 
\item For some results we will refer to the condition
\begin{equation}
\label{eq:Tp0ser}
T_{p0} < \infty,
\end{equation}
where $T_{p0} \triangleq \Ex_p[\tau_0] = \sum_{n=1}^\infty p_n T_{n0}$ is
the \emph{mean first-passage time to $0$ starting from steady state}.
Note that from the identities for $T_{n0}$ in (\ref{eq:Tn0}) we have
\begin{equation}
\label{eq:Tp0series}
T_{p0}  = \sum_{n=0}^\infty \frac{\bar{P}_n^2}{\lambda_n p_n}.
\end{equation}
 See also \cite[Eq.\ (6.6)]{coolenVanDoorn02} and  \cite[Theorem 4 and Eq.\ (3.6)]{karlMcG57}. 
Furthermore, 
\cite[Theorem 6.1]{coolenVanDoorn02} shows that satisfaction of (\ref{eq:Tp0ser})  is equivalent to existence of the \emph{deviation matrix} $D = (d_{mn})$  of the Markov chain, defined by
\[
d_{mn} \triangleq \int_0^\infty (p_{mn}(t) - p_n) \, dt,
\] 
where $p_{mn}(t) \triangleq \mathbb{P}_m\{X(t) = n\}$ and $\mathbb{P}_m$ is the probability starting from $m$. 
Note that the bias is given in terms of $D$ by $\beta_m = \sum_{n} d_{mn} c_n$, i.e., 
$\beta = D c$.
\item The quantities $T_{\infty 0}$ and $T_{p0}$ are related by 
\begin{equation}
\label{eq:Tinf0Tp0}
T_{\infty 0} = \sum_{n=0}^\infty \bar{P}_n T_n^+ + T_{p0},
\end{equation}
which readily follows from (\ref{eq:tinfty0}), (\ref{eq:Tp0series}), (\ref{eq:Tn}) and $P_n - \bar{P}_n = P_n^2 - \bar{P}_n^2$.
In light of (\ref{eq:Tinf0Tp0}), we define 
\begin{equation}
\label{eq:Tinfp}
T_{\infty p} \triangleq 
\sum_{n=0}^\infty \bar{P}_n T_n^+,
\end{equation}
 which represents the \emph{mean first passage-time to steady state starting at $\infty$}.
\end{enumerate}
\end{remark}

The following result shows that the ratios $T_{n0}/T_{0n}$ strictly decrease to $0$.
\begin{lemma}
\label{lma:ratioTn00n} \textup{ }
\begin{enumerate}[(a)]
\item $\displaystyle \frac{\bar{P}_n}{P_n} < \frac{T_{n+1, 0}}{T_{0, n+1}} < \frac{T_{n 0}}{T_{0n}},$ for $n \in \mathbb{N};$
\item $\displaystyle \frac{T_{n 0}}{T_{0n}} \searrow 0$ as $n \to \infty.$
\end{enumerate}
\end{lemma}
\begin{proof} (a)
Using (\ref{eq:T0np1H0np1}), (\ref{eq:Tn0}) and the mediant inequality (see (\ref{eq:ein1})) we obtain 
\begin{equation}
\label{eq:ratbTnTn}
\frac{\bar{P}_n}{P_n} < \frac{T_{n+1, 0}}{T_{0, n+1}} < \frac{T_{n 0}}{T_{0n}} \; \textup{ if and only if } \; \frac{\bar{P}_n}{P_n} < \frac{T_{n 0}}{T_{0n}}.
\end{equation}
We next show by induction that $\bar{P}_n/P_n < T_{n+1, 0}/T_{0, n+1}$ for $n \geqslant 1$. For $n = 1$, 
\[
\frac{T_{10}}{T_{01}} = \frac{\bar{P}_0}{P_0} > \frac{\bar{P}_1}{P_1},
\] 
whence the result follows by (\ref{eq:ratbTnTn}).
Suppose now the result holds for $n-1$. Then, 
\begin{equation}
\label{eq:fortn00n}
\frac{T_{n 0}}{T_{0n}} > \frac{\bar{P}_{n-1}}{P_{n-1}} > \frac{\bar{P}_{n}}{P_{n}},
\end{equation}
which, using (\ref{eq:ratbTnTn}), gives the result for $n$. This completes the induction.
Now, part (a) follows since, by (\ref{eq:ratbTnTn}), it is a consequence of (\ref{eq:fortn00n}).

(b) 
Since
\[
\frac{T_{n 0}}{T_{0n}} = \frac{\sum_{j =0}^{n-1} \frac{\bar{P}_j}{\lambda_j p_j}}{\sum_{j =0}^{n-1} \frac{P_j}{\lambda_j p_j}},
\]
with $T_{0n}$ strictly increasing and divergent and $\bar{P}_{n-1}/P_{n-1}  \to 0$, 
the \emph{Stolz--Ces\`aro theorem} (see \cite[Theorem 1.22]{muresan08}) gives that $T_{n 0}/T_{0n} \to 0$.
\end{proof}
 
By analogy with expression (\ref{eq:Zirel}) for $Z_n$, we further define the cost ratios
\begin{equation}
\label{eq:Z0np1}
Z_{0n} \triangleq \frac{H_{0n}}{T_{0n}} \quad
\textup{and}
\quad
Z_{n0} \triangleq \frac{H_{n0}}{T_{n0}}, \quad n \in \mathbb{N}.
\end{equation}
Since $Z_n \to \zeta$ as $n \to \infty$, this raises the question of whether $Z_{0n}$ and $Z_{n0}$ converge to $\zeta$.
We will settle this in the affirmative for $Z_{0n}$.
We need a preliminary result on invariance relations for certain ratios of mean costs to mean times.
\begin{lemma}
\label{lma:ratioHpmTpm} \textup{ }
\begin{itemize}
\item[\textup{(a)}] $\displaystyle \frac{H_n^+ + H_{n+1}^-}{T_n^+ + T_{n+1}^-} = \zeta,$ for $n \in \mathbb{N}_0;$
\item[\textup{(b)}] $\displaystyle 
\frac{H_{0n} + H_{n0}}{T_{0n} + T_{n0}} = \zeta,$ for $n \in \mathbb{N}.$
\end{itemize}
\end{lemma}
\begin{proof}
(a) Using (\ref{eq:Tn}), (\ref{eq:Hn}) and Lemma \ref{lma:CPm}, we have
\[
\frac{H_n^+ + H_{n+1}^-}{T_n^+ + T_{n+1}^-} = 
\frac{\frac{C_n}{\lambda_n p_n} + \frac{\bar{C}_n}{\lambda_n p_n}}{\frac{P_n}{\lambda_n p_n} + \frac{\bar{P}_n}{\lambda_n p_n}} = 
\zeta.
\]

(b) Using (\ref{eq:T0np1H0np1}) and (\ref{eq:Tn0}) we can write
\[
\frac{H_{0n} + H_{n0}}{T_{0n} + T_{n0}} = 
\frac{\sum_{j=0}^{n-1} \frac{C_j}{\lambda_j p_j} + \sum_{j=0}^{n-1} \frac{\bar{C}_j}{\lambda_j p_j}}{\sum_{j=0}^{n-1} \frac{P_j}{\lambda_j p_j} + \sum_{j=0}^{n-1} \frac{\bar{P}_j}{\lambda_j p_j}} = 
\frac{\sum_{j=0}^{n-1} \frac{\zeta}{\lambda_j p_j}}{\sum_{j=0}^{n-1} \frac{1}{\lambda_j p_j}} = \zeta.
\]
\end{proof}

We can now prove the following result on the convergence of $Z_{0n}$.

\begin{lemma}
\label{lma:Z0np1}
$Z_{0n} \to \zeta$ as $n \to \infty.$ 
\end{lemma}
\begin{proof}
Since
\[
Z_{0n} = \frac{H_{0n}}{T_{0n}} = 
\frac{H_0^+ + \cdots + H_{n-1}^+}{T_0^+ + \cdots + T_{n-1}^+},
\]
with $T_{0n}$ strictly increasing and divergent and $H_{n-1}^+/T_{n-1}^+ = Z_n \to \zeta$, 
the Stolz--Ces\`aro theorem (see \cite[Theorem 1.22]{muresan08}) gives that $Z_{0n} \to \zeta$.
\end{proof} 

\begin{remark}
\label{re:Zn0m} 
From  Lemma \ref{lma:ratioHpmTpm}(b) and (\ref{eq:ei1}) we can write
\[
\zeta = \frac{H_{n0} + H_{0n}}{T_{n0} + T_{0n}} = 
Z_{n0} + \frac{T_{0n}}{T_{n0} + T_{0n}} (Z_{0n} - Z_{n0}),
\]
hence
\begin{equation}
\label{eq:Zn0ref}
Z_{n0} = \zeta -  \frac{T_{0n}}{T_{n0}} (Z_{0n} - \zeta).
\end{equation}
Thus, the asymptotic behavior of $Z_{n0}$ depends on that of the scaled tail of $Z_{0n}$ in (\ref{eq:Zn0ref}). See Lemmas \ref{lma:ratioTn00n}(b) and \ref{lma:Z0np1}. Typically, $Z_{n0}$ will \emph{not} converge to $\zeta$. 
\end{remark}

\section{Exact solution: explicit expressions and properties}
\label{s:ersp}
\subsection{Explicit expressions for the solution to Poisson's equation}
\label{s:estpe}
We next derive explicit expressions for the exact solution to Poisson's equation (\ref{eq:pe1}).
We start by obtaining expressions formulated 
 in terms of the mean \emph{upward} first-passage times and costs considered in \S\ref{s:bdma}.

Let $b(z; a)$ be the solution to (\ref{eq:pe1}) for given $z$ when $b_0 = a$.
To analyze $b(z; a)$, 
consider the reformulation (\ref{eq:pe2}) of (\ref{eq:pe1}). Note that
\begin{equation}
\label{eq:recbj}
b_{n} = b_{0} + \sum_{j=0}^{n-1} f_{j}, \quad n \in \mathbb{N}. 
\end{equation}
Let $f(z)$ be
the solution to (\ref{eq:pe2}), which is constructed by the forward recurrence
\begin{equation}
\label{eq:recsolpe2}
f_0(z) = 
\frac{z - c_0}{\lambda_0}, \quad
f_n(z) = 
\frac{z - c_n}{\lambda_n} + \frac{\mu_n}{\lambda_n} f_{n-1}(z), \quad n = 1, 2, \ldots
\end{equation}

The following result  represents $f_n(z)$ and $b_n(z; a)$ as affine functions of $z$.

\begin{proposition}
\label{pro:fzexp} For $n \in \mathbb{N}_0,$
\begin{itemize}
\item[\textup{(a)}]
$\displaystyle f_n(z) =  T_n^+ z  - H_n^+ =  T_n^+ (z  - Z_n);$
\item[\textup{(b)}]
$\displaystyle b_n(z; a) =   a + T_{0n} z  - H_{0n} = a + T_{0n} (z  - Z_{0n}).$
\end{itemize}
\end{proposition}
\begin{proof}
(a)
Since $T^+$ and $H^+$ satisfy (\ref{eq:Pi}) and (\ref{eq:Hi}),  $f = T^+ z - H^+$ satisfies (\ref{eq:recsolpe2}), whence $f(z) = T^+ z - H^+ = T^+ (z - Z)$, where we have also used (\ref{eq:Zirel}).

(b) We have $b_{0}(z; a) = a$. For $n \geqslant 1$,  from (\ref{eq:recbj}), part (a),  (\ref{eq:T0np1H0np1}) and (\ref{eq:Z0np1}),
\[
b_{n}(z; a) = a + \sum_{j=0}^{n-1} f_{j}(z) = 
a + \sum_{j=0}^{n-1} (T_j^+ z  - H_j^+) = 
a + T_{0n} z  - H_{0n} = a + T_{0n} (z  - Z_{0n}).
\]
\end{proof}  

From  Proposition \ref{pro:fzexp} we immediately obtain the following result, which gives expressions for the $\varphi_n = f_n(\zeta)$ and $b_n = b_n(\zeta; a)$ solving  (\ref{eq:varphirec}) and (\ref{eq:pe0}).

\begin{corollary}
\label{cor:phiexp}
For $n \in \mathbb{N}_0,$ 
\begin{itemize}
\item[\textup{(a)}]
$\displaystyle \varphi_n =   T_n^+ (\zeta  - Z_n);$
\item[\textup{(b)}]
$\displaystyle b_n =   b_0 + T_{0n} (\zeta  - Z_{0n}).$
\end{itemize}
\end{corollary}

\begin{remark}
\label{re:phiexp}
\begin{itemize}
\item[\textup{(a)}]
The identity in Corollary \ref{cor:phiexp}(a) was given in \cite[Eq.\ (6)]{krish90} for a particular queueing model, the M/M/$m$ queue with $c_n = n$.
\item[\textup{(b)}] Corollary \ref{cor:phiexp}(b) is consistent with known results for discrete-time Markov chains. See, e.g., \cite[Prop. A.3.1(ii)]{meyn07} and \cite[Theorem 1]{tijms17}.
\end{itemize}
\end{remark}

The following result gives alternate expressions for the exact solution to Poisson's equations (\ref{eq:varphirec}) and (\ref{eq:pe0}), which are formulated in terms of the mean \emph{downward} first-passage times and costs considered in \S\ref{s:bdma}.

\begin{proposition}
\label{pro:bfzexp} \textup{ }
\begin{itemize}
\item[\textup{(a)}]
$\displaystyle \varphi_n =  T_{n+1}^-  (Z_{n+1}^- - \zeta),$ for $n \in \mathbb{N}_0;$
\item[\textup{(b)}]
$\displaystyle b_n =   b_0 + T_{n0} (Z_{n0} - \zeta),$ for $n \in \mathbb{N}.$
\end{itemize}
\end{proposition}
\begin{proof}
(a)
From Corollary \ref{cor:phiexp}(a), (\ref{eq:Zirel}), Lemma \ref{lma:ratioHpmTpm}(a) and (\ref{eq:Zirelm}), we obtain 
\[
\varphi_n =   T_n^+ (\zeta  - Z_n) = T_n^+ \zeta  - H_n^+ = H_{n+1}^- - \zeta \, T_{n+1}^- = T_{n+1}^- (Z_{n+1}^- - \zeta).
\]

(b) From Corollary \ref{cor:phiexp}(b), (\ref{eq:Z0np1}), and Lemma \ref{lma:ratioHpmTpm}(b), we obtain 
\[
b_n -  b_0 = T_{0n} (\zeta  - Z_{0n}) = \zeta \, T_{0n}   - H_{0n} = H_{n0} - \zeta \, T_{n0} = T_{n0} (Z_{n0} - \zeta).
\]
\end{proof}  

\subsection{Convexity of relative cost function}
\label{s:pofzi}
We next draw on the exact expressions derived above to obtain the practically relevant structural result of solutions to Poisson's equation (\ref{eq:pe0}) stated in Theorem \ref{the:pofzi}.
The short proof given next  draws on substantial preliminary groundwork, which is laid in Appendix \ref{a:pofzi}.

\begin{proof}[Proof of Theorem \ref{the:pofzi}]
From Corollary \ref{cor:phiexp}(a), we immediately obtain 
\[
\Delta \varphi_n = \zeta \Delta T_n^+ - \Delta H_n^+,
\]
and hence, since $\Delta T_n^+ > 0$ by Lemma \ref{lma:PHstarinc},
\begin{equation}
\label{eq:varphind}
\Delta \varphi_n \geqslant 0 \quad \textup{if and only if} \quad \frac{\Delta H_n^+}{\Delta T_n^+} \leqslant \zeta.  
\end{equation}

Now, Lemmas \ref{lma:Deltaziinc} and \ref{lma:limDeltaziinc} ensure that $\Delta H_n^+/\Delta T_n^+$ grows to $\zeta$ as $n \to \infty$. In light of (\ref{eq:varphind}), this completes the proof.
\end{proof}

\subsection{Bias and asymptotic variance: exact expressions in terms of $\varphi$}
\label{s:ebav}
We next turn to exact evaluation of the bias $\beta$ and the asymptotic variance $\sigma^2$. 
The following result gives new expressions 
in terms of the marginal relative cost $\varphi$.
We assume that $\beta$ is characterized by (\ref{eq:adc1}), and that $\sigma^2$ is well defined and finite, being given by
(\ref{eq:sigma22}). To ensure the validity of interchanging the order of summation in certain series arising in the proofs, we further assume that cost rates $c_n$ are nonnegative and nondecreasing, so that Assumption \textup{\ref{ass:hmulambda}(ii.a)} holds.

\begin{proposition}
\label{pro:betasigma2phi} 
Suppose that $c$ satisfies Assumption \textup{\ref{ass:hmulambda}(ii.a)}$.$    Then
\begin{itemize}
\item[\textup{(a)}] $\displaystyle \beta_{0} = - \sum_{j =0}^{\infty} \bar{P}_j \varphi_j, \quad \beta_{n} = \beta_{0} +\sum_{j=0}^{n-1} \varphi_j,
 \quad n \in \mathbb{N};$
\item[\textup{(b)}] $\displaystyle \sigma^2 = 2 \sum_{n=0}^{\infty} \lambda_n   p_n \varphi^2_n.$
\end{itemize}
\end{proposition}
\begin{proof}
(a) From (\ref{eq:recbj}) we have
$\beta_{n} = \beta_{0} + \sum_{j=0}^{n-1} \varphi_j$ and, using (\ref{eq:adc1}), we can write
\[
0 = \sum_{n =0}^{\infty} \beta_n p_n = \sum_{n =0}^{\infty} (\beta_{0} + \sum_{j=0}^{n-1} \varphi_j) p_n =
\beta_{0} + \sum_{j=0}^{\infty} \varphi_j \sum_{n =j+1}^{\infty}  p_n = \beta_{0} + \sum_{j=0}^{\infty} \bar{P}_j \varphi_j,
\]
where the interchange in the order of summation is justified by Tonelli's theorem, since $\varphi \geqslant 0$ by Lemma \ref{lma:ziinc}(b). 
Therefore, 
\begin{equation}
\label{eq:beta0eq}
\beta_{0} = - \sum_{j =0}^{\infty} \bar{P}_j \varphi_j.
\end{equation}

(b) We can write
\begin{align*}
\sigma^2 & = 
2 \sum_{n =0}^\infty \beta_{n} c_{n} p_{n}  = 2 \sum_{n =0}^\infty \bigg(\beta_{0} + \sum_{j=0}^{n-1} \varphi_j\bigg) c_{n} p_{n}  \\
& = 2  \beta_{0} \zeta + 2 \sum_{j=0}^{\infty} \varphi_j\sum_{n=j+1}^{\infty}  c_{n} p_{n}  = 2  \beta_{0} \zeta  + 2 \sum_{j=0}^{\infty} \varphi_j \bar{C}_j \\
&  = 
2 \sum_{n=0}^{\infty} (\bar{C}_n - \zeta \bar{P}_n) \varphi_n = 2 \sum_{n=0}^{\infty} \lambda_n p_n (H_{n+1}^- - T_{n+1}^-) \varphi_n  = 2 \sum_{n=0}^{\infty} \lambda_n p_n \varphi^2_n,
 \end{align*}
using in turn (\ref{eq:sigma22}), part (a),  (\ref{eq:PCZ1is}), Lemma \ref{lma:CPm} and Proposition \ref{pro:bfzexp}(a), and the interchange in the order of summation is justified as in part (a).
\end{proof}

\section{Approximate numerical solution by forward recurrence}
\label{s:fraeas}
This section presents an error analysis of standard forward recurrence for the numerical solution of Poisson's equation.
It further addresses how the resulting approximation errors in the computed solution affect the accuracy of computed approximations to the bias $\beta$ and the asymptotic variance $\sigma^2$. 

Given an approximation $\widehat{x}$ to a number $x$, 
we denote by $E_{\textup{abs}}(\widehat{x}) \triangleq \widehat{x} - x$ and 
$E_{\textup{rel}}(\widehat{x}) \triangleq E_{\textup{abs}}(\widehat{x})/x$ the corresponding 
approximation errors in absolute and relative terms, respectively, provided that $x \neq 0$ in the latter case. 

Note that, typically, $\widehat{x}$ will be a floating-point approximation to $x$, and hence its relative error will be bounded as (see, e.g., \cite[Ch.\ 2]{higham02})
\begin{equation}
\label{eq:ieeesa}
|E_{\textup{rel}}(\widehat{x})| < u,
\end{equation}
where $u$ is the \emph{unit roundoff}. Thus, in IEEE standard arithmetic, $u = 2^{-24} \approx 5.96 \times 10^{-8}$ for single precision, and $u = 2^{-53} \approx 1.11 \times 10^{-16}$ for double precision. 

\subsection{Approximate numerical evaluation of $\varphi$ and $b$: Error amplificaton factors}
\label{s:nephib}
We start by addressing the approximate numerical evaluation of   relative costs $b_n$ and marginal relative costs $\varphi_n$, using an approximate input $z = \widehat{\zeta} \neq  \zeta$ instead of $z = \zeta$ in (\ref{eq:pe1}) and (\ref{eq:pe2}), respectively.
We will assume that to be the only source of error, so the computed approximation to $\varphi$ is
$\widehat{\varphi} = f(\widehat{\zeta})$. See (\ref{eq:recsolpe2}).

We next address how the approximation errors in the input $\widehat{\zeta}$ are amplified to  corresponding errors in the computed outputs $\widehat{\varphi}_n$ and $\widehat{b}_n$, where the latter are computed from $\widehat{\varphi}$ through (\ref{eq:recbj}) with $\widehat{b}_0 = b_0 = a$.
The following result  identifies the corresponding \emph{error amplificaton factors}.

Note that, as is standard in error analysis of numerical algorithms (see, e.g., \cite{gautschi61}), we call $A_n$ the
absolute (resp.\ relative) \emph{error amplification factor} of a computed quantity, such as 
$\widehat{\varphi}_n$, with respect to the error in the approximate input, which is $\widehat{\zeta}$ in our case, if  $E_{\textup{abs}}(\widehat{\varphi}_{n}) =  A_n E_{\textup{abs}}(\widehat{\zeta})$ (resp.\ $E_{\textup{rel}}(\widehat{\varphi}_{n}) =  A_n E_{\textup{rel}}(\widehat{\zeta})$).

We assume henceforth that relative errors  are well defined, i.e., $\zeta, \varphi_n, b_n \neq 0$.
 
\begin{proposition}[Error amplification factors]
\label{pro:phihatsol}
\textup{ }
\begin{itemize}
\item[\textup{(a)}]
$\displaystyle E_{\textup{abs}}(\widehat{\varphi}_{n}) =  T_n^+ E_{\textup{abs}}(\widehat{\zeta}),$ for $n \in \mathbb{N}_0;$
\item[\textup{(b)}]
$\displaystyle E_{\textup{rel}}(\widehat{\varphi}_{n})  = \zeta \frac{T_n^+}{\varphi_n} E_{\textup{rel}}(\widehat{\zeta}) =  \frac{\zeta}{\zeta  - Z_n} E_{\textup{rel}}(\widehat{\zeta}),$ for $n \in \mathbb{N}_0;$
\item[\textup{(c)}]
$\displaystyle E_{\textup{abs}}(\widehat{b}_{n}) =
T_{0n} \, E_{\textup{abs}}(\widehat{\zeta}),$ for $n \in \mathbb{N};$
\item[\textup{(d)}]
$\displaystyle E_{\textup{rel}}(\widehat{b}_{n}) =
\zeta \frac{T_{0n}}{b_n} \, E_{\textup{rel}}(\widehat{\zeta}) =
\frac{\zeta}{{a}/{T_{0n}} + \zeta - Z_{0n}} \, E_{\textup{rel}}(\widehat{\zeta}),$ for $n \in \mathbb{N}.$
  
\end{itemize}
\end{proposition}
\begin{proof}
All parts follow straightforwardly from Proposition \ref{pro:fzexp}. Thus, e.g., for part (a), taking $z = \widehat{\zeta}$ and $z = \zeta$ in Proposition \ref{pro:fzexp}(a) gives 
$\widehat{\varphi}_n = T_n^+ (\widehat{\zeta} - Z_n)$ and $\varphi_n = T_n^+ (\zeta - Z_n)$, hence $\widehat{\varphi}_n - \varphi_n = T_n^+ (\widehat{\zeta} -\zeta)$.
\end{proof}
  
The following result elucidates the numerical instability phenomenon (cf.\ Table \ref{tab:numex1} in \S\ref{s:intro}), by clarifying the asymptotic behavior of the error amplification factors in Proposition \ref{pro:phihatsol}. Note that $\sgn(x) \in \{-1, 0, 1\}$ denotes the sign of $x$.

\begin{corollary}
\label{cor:phihatsol}
As $n \to \infty,$
\begin{itemize}
\item[\textup{(a.1)}] 
$E_{\textup{abs}}(\widehat{\varphi}_{n}) \to \sgn(E_{\textup{abs}}(\widehat{\zeta})) \cdot \infty$ \enspace if and only if \enspace
$\lambda_n p_n \to 0;$ 
\item[\textup{(a.2)}] 
$\sup_n E_{\textup{abs}}(\widehat{\varphi}_{n}) = \infty$ \enspace if and only if \enspace
$\inf_n \lambda_n p_n = 0;$ 
\item[\textup{(b)}]
$\displaystyle \left|E_{\textup{rel}}(\widehat{\varphi}_{n})\right| \to  \infty;$
\item[\textup{(c)}]
$\displaystyle E_{\textup{abs}}(\widehat{b}_{n}) \to \sgn(E_{\textup{abs}}(\widehat{\zeta}))  \cdot \infty;$
\item[\textup{(d)}]
$\displaystyle |E_{\textup{rel}}(\widehat{b}_{n})| \to  \infty.$
\end{itemize}
\end{corollary}
\begin{proof}
The results follow from (\ref{eq:Tncons}), (\ref{eq:supTncons}), 
 Proposition \ref{pro:phihatsol}, Lemma \ref{lma:tnpinfty}, (\ref{eq:PCZ1is}), Remark \ref{re:T0np1div}(a) and Lemma \ref{lma:Z0np1}.
 \end{proof}
 
 \begin{remark}
\label{re:phihatsol} 
\hspace{1in}
\begin{itemize}
\item[\textup{(a)}] Corollary \ref{cor:phihatsol}(b, d) shows that  forward recurrence is inherently numerically unstable, in that the magnitudes of \emph{relative errors} in the computed approximations $\widehat{\varphi}_{n}$ and $\widehat{b}_{n}$ to $\varphi_{n}$ and $b_{n}$ diverge to infinity as $n \to \infty$, as they are inversely proportional to the asymptotically vanishing \emph{tails} $\zeta - Z_n$ and $\zeta - Z_{n0}$ (approximately so for $\widehat{b}_{n}$ when $a \neq 0$), respectively. Corollary \ref{cor:phihatsol}(c) shows that such is also the case for the approximation error of $\widehat{b}_{n}$.
\item[\textup{(b)}] It is of interest to consider how the relative approximation errors of $\widehat{\varphi}_{n}$ and $\widehat{b}_{n}$ grow asymptotically in relation to $1/p_n$. From Proposition \ref{pro:phihatsol}(b, d) it follows that $p_n E_{\textup{rel}}(\widehat{\varphi}_{n})$ and $p_n E_{\textup{rel}}(\widehat{b}_{n})$ are asymptotically proportional to $p_n/(\zeta  - Z_n)$ and $p_n/(\zeta - Z_{0n})$, respectively. Note that the latter ratios might, e.g., converge to a finite limit or diverge to infinity. 
\item[\textup{(c)}] The approximation error of $\widehat{\varphi}_{n}$ may be bounded or unbounded. 
Corollary \ref{cor:phihatsol}(a.2) shows that it is bounded if and only if $\lambda_n p_n$ is bounded away from $0$. As an example where $E_{\textup{abs}}(\widehat{\varphi}_{n})$ is bounded, take $\lambda_n \triangleq \lambda \rho^{-n-1}$ and $\mu_n \triangleq \mu \rho^{-n-1}$, with $\rho \triangleq \lambda/\mu < 1$. Then, $p_n = (1-\rho) \rho^n$, so $\lambda_n p_n \equiv   \mu-\lambda$. 
\item[\textup{(d)}] If arrival rates $\lambda_n$ are bounded above, so that $\bar{\lambda} \triangleq \sup_n \lambda_n < \infty$, as is often the case in applied models, Proposition \ref{pro:phihatsol}(a) and (\ref{eq:Pi}) give that, for large $n$, 
$|E_{\textup{abs}}(\widehat{\varphi}_{n})| \approx   |E_{\textup{abs}}(\widehat{\zeta})| / (\lambda_n p_n) \geqslant |E_{\textup{abs}}(\widehat{\zeta})| / (\bar{\lambda} p_n) $. In such a case, $\widehat{\varphi}_{n}$ will substantially deviate from $\varphi_{n}$ for unlikely states since, asymptotically, $|E_{\textup{abs}}(\widehat{\varphi}_{n})|$ will grow at least inversely proportional to $p_n$.
\item[\textup{(e)}] Corollary \ref{cor:phihatsol}(a.1) explains the opposite signs in the behavior for large $n$ of the computed approximations $\widehat{\varphi}_{n}$ and $\widetilde{\varphi}_{n}$ in  Table \ref{tab:numex1} (see \S \ref{s:intro}).
\end{itemize}
\end{remark} 

The following result highlights the strikingly different asymptotic behavior of
the computed  $\widehat{\varphi}_{n}$ and $\widehat{b}_{n}$ compared to the exact $\varphi_{n}$ and $b_{n}$.
Note that the notation  $y_{n} = \Theta(x_n)$ means that $y_{n}$ is \emph{asymptotically proportional} to $x_n$.
 
\begin{corollary}
\label{cor:2phihatsol}
Let $\widehat{\zeta} \neq \zeta \neq 0.$
As $n \to \infty,$
\begin{itemize}
\item[\textup{(a)}] If $\lambda_n p_n \to 0,$ $\widehat{\varphi}_{n} = \Theta({T_n^+})$ and $\varphi_{n} = o({T_n^+}),$ and hence $\varphi_{n} = o(\widehat{\varphi}_{n});$
\item[\textup{(b)}]
$\widehat{b}_{n} = \Theta({T_{0n}})$ and $b_{n} = o({T_{0n}}),$ and hence $b_{n} = o(\widehat{b}_{n}).$ 
\end{itemize}
\end{corollary}
\begin{proof}
(a) By Proposition \ref{pro:fzexp}(a), 
$\widehat{\varphi}_{n} = T_n^+ (\widehat{\zeta}-Z_n)$ and
$\varphi_{n} = T_n^+ (\zeta-Z_n)$, which gives the result using that $Z_n \to \zeta$ and $T_n^+ \to \infty$ (see  (\ref{eq:PCZ1is}) and Lemma \ref{lma:tnpinfty}), since 
$\widehat{\varphi}_{n}/T_n^+ \to \widehat{\zeta}-\zeta \neq 0$ and 
$\varphi_{n}/T_n^+ \to 0$ as $n \to \infty$.

(b)
The result follows similarly as part (a) using Proposition \ref{pro:fzexp}(b) and that $Z_{0, n+1} \to \zeta$ and $T_{0, n+1} \to \infty$ (see  Lemma \ref{lma:Z0np1} and Remark \ref{re:T0np1div}(a)).
\end{proof}

\subsection{Error analysis of bias and asymptotic variance computation}
\label{s:eacrvav}
We next address how the approximation errors in the computed $\widehat{\varphi}_n$ resulting from errors in input $\widehat{\zeta}$, as characterized in Proposition \ref{pro:phihatsol}, affect the accuracy of the 
computed bias $\widehat{\beta}$ and asymptotic variance $\widehat{\sigma}^2$. 

Specifically, we discuss next whether the expressions given in Proposition \ref{pro:betasigma2phi} for $\beta_{n}$ and  $\sigma^2$ are also valid 
to approximately evaluate such quantities, substituting $\widehat{\varphi}$ for 
$\varphi$. 
We would thus 
obtain approximations $\widehat{\beta}$ and $\widehat{\sigma}^2$ given by
\begin{equation}
\label{eq:approxbs}
\begin{split}
 \widehat{\beta}_{0} & = - \sum_{j =0}^{\infty} \bar{P}_j \widehat{\varphi}_{j}, \quad 
 \widehat{\beta}_{n}  = \widehat{\beta}_{0} + \sum_{j=0}^{n-1} \widehat{\varphi}_{j}, \quad n \in \mathbb{N}, \\
\widehat{\sigma}^2 & = 2 \sum_{n=0}^{\infty} \lambda_n p_n \widehat{\varphi}^2_n.
\end{split}
\end{equation}

Yet, Proposition \ref{pro:nbetasigma2phi} below shows that the  expressions in (\ref{eq:approxbs}) will typically be invalid,  due to divergence of the infinite series involved. 
We need a preliminary result.
See Remark \ref{re:T0np1div} and the expressions for $T_{\infty 0}$ and $T_{\infty p}$ in (\ref{eq:tinfty0}) and (\ref{eq:Tinfp}).

 \begin{lemma}
\label{lma:prelr} 
\hspace{1in}
\begin{itemize}
\item[\textup{(a)}] $\displaystyle \sum_{n =0}^{\infty} P_n T_n^+ = \infty;$
\item[\textup{(b)}] $\displaystyle T_{\infty p} = \infty$ \enspace if and only if \enspace $T_{\infty 0} = \infty;$
\item[\textup{(c)}]
$ \displaystyle
\frac{\sum_{n =0}^{N} P_n \varphi_n}{\sum_{n =0}^{N} P_n T_n^+} \to 0 \enspace \textup{as} \enspace N \to \infty.
$
\end{itemize}
\end{lemma}
\begin{proof}
(a) This part follows from $\sum_{n =0}^{\infty} T_n^+ = \infty$ (see Remark \ref{re:T0np1div}(a)) and the limit comparison test, using that  $P_n \to 1$ as $n \to \infty$.

(b) By the limit comparison test, $T_{\infty p} = \sum_{n =0}^{\infty} \bar{P}_n T_n^+ = \infty$ if and only if $T_{\infty 0} = \sum_{n =0}^{\infty} \bar{P}_n /\lambda_n p_n = \infty$, using that, by (\ref{eq:Tn}),  $\lambda_n p_n T_n^+ = P_n \to 1$ as $n \to \infty$.

(c)
This part follows from the Stolz--Ces\`aro theorem (see \cite[Theorem 1.22]{muresan08}), since $\sum_{n =0}^{N} P_n T_n^+ \nearrow \infty$ strictly and 
$\varphi_n/T_n^+ = \zeta - Z_n \to 0$  as $n \to \infty$.
\end{proof}

In practice, one would approximate the infinite series in Proposition \ref{pro:betasigma2phi} by truncation. Thus, 
let
$\beta_{0, N} \triangleq 
-\sum_{n =0}^{N} \bar{P}_n \varphi_{n}$ and 
$\sigma_N^2 \triangleq  2 \sum_{n =0}^{N} \lambda_n p_n \varphi^2_n$,  which 
converge to $\beta_0$ and $\sigma^2$ as $N \to \infty$. Consider also the computed approximations
$\widehat{\beta}_{0, N} \triangleq 
-\sum_{n =0}^{N} \bar{P}_n \widehat{\varphi}_{n}$ and
$\widehat{\sigma}_N^2 \triangleq  2 \sum_{n =0}^{N} \lambda_n p_n \widehat{\varphi}^2_n$.
The following result gives expressions for the latter quantities and elucidates their asymptotic behavior.

\begin{proposition}
\label{pro:nbetasigma2phi} For $N \in \mathbb{N}_0,$ 
\begin{itemize}
\item[\textup{(a)}] $\displaystyle \widehat{\beta}_{0, N} = \beta_{0, N} - E_{\textup{abs}}(\widehat{\zeta}) \sum_{n =0}^{N} \bar{P}_n T_n^+;$
\item[\textup{(b)}] if $T_{\infty 0} = \infty$ then $|\widehat{\beta}_{0, N}| \to \infty;$
otherwise$,$ 
$\displaystyle \widehat{\beta}_{0, N} \to \beta_{0} - T_{\infty p} E_{\textup{abs}}(\widehat{\zeta});$ 
\item[\textup{(c)}] 
$\displaystyle \widehat{\sigma}_N^2 = \sigma_N^2  + 
 4 E_{\textup{abs}}(\widehat{\zeta}) \sum_{n =0}^{N} P_n \varphi_n 
 + 2
 E_{\textup{abs}}^2(\widehat{\zeta}) \sum_{n =0}^{N} P_n T_n^+;$
\item[\textup{(d)}] 
 $\displaystyle  \widehat{\sigma}_N^2 \to \infty$ as $N \to \infty.$
\end{itemize}
\end{proposition}
\begin{proof}
(a) From Proposition \ref{pro:phihatsol}(a), i.e., 
 $\widehat{\varphi}_{n} = \varphi_n + T_n^+ E_{\textup{abs}}(\widehat{\zeta})$, we have
 \[
\widehat{\beta}_{0, N} = -\sum_{n =0}^{N} \bar{P}_n \widehat{\varphi}_{n} = 
 -\sum_{n =0}^{N} \bar{P}_n \varphi_n - E_{\textup{abs}}(\widehat{\zeta}) \sum_{n =0}^{N} \bar{P}_n T_n^+ = \beta_{0, N} -E_{\textup{abs}}(\widehat{\zeta}) \sum_{n =0}^{N} \bar{P}_n T_n^+.
\]

(b) The result follows by letting $N \to \infty$ in part (a), using Lemma \ref{lma:prelr}(b).

(c)
Using Proposition \ref{pro:phihatsol}(a)  and 
 (\ref{eq:Tn}), we obtain
\begin{equation}
\label{eq:ninfvar}
\begin{split}
\widehat{\sigma}_N^2 & =  2 \sum_{n =0}^{N} \lambda_n p_n \widehat{\varphi}^2_n =  2 \sum_{n =0}^{N} \lambda_n p_n (\varphi_n +  E_{\textup{abs}}(\widehat{\zeta}) T_n^+)^2 
\\ & = 
 \sigma_N^2 + 
 4 E_{\textup{abs}}(\widehat{\zeta}) \sum_{n =0}^{N} \lambda_n p_n T_n^+ \varphi_n  +
 2 E_{\textup{abs}}^2(\widehat{\zeta}) \sum_{n =0}^{N} \lambda_n p_n (T_n^+)^2 \\
& = \sigma_N^2  + 
 4 E_{\textup{abs}}(\widehat{\zeta}) \sum_{n =0}^{N} P_n \varphi_n 
 +
 2 E_{\textup{abs}}^2(\widehat{\zeta}) \sum_{n =0}^{N} P_n T_n^+. 
\end{split}
\end{equation}

(d) The result follows from part (c) and Lemma \ref{lma:prelr}(c).
\end{proof}

\section{Approximate numerical solution by backward recurrence}
\label{s:abr}
Instead of computing the approximate solution to Poisson's equation 
through standard forward recurrence, we now consider using \emph{backward recurrence}, which to the author's knowledge has not been previously explored for this model.

This approach is based on the observation that the $f_n(z)$ in \S\ref{s:estpe} satisfy the backward recurrence (cf.\ (\ref{eq:recsolpe2}))\begin{equation}
\label{eq:fnm1zbr}
f_{n-1}(z) = \frac{c_n - z}{\mu_n} + \frac{\lambda_n}{\mu_n} f_n(z), \quad n \in \mathbb{N}.
\end{equation}

Given an approximation $\widehat{\zeta} \neq \zeta$, 
 fix a large integer $N$ and set $\widehat{\varphi}^{N}_{N}$ to an approximation to $\varphi^{\mathstrut}_{N}$, which we require to satisfy the asymptotic condition  
\begin{equation}
\label{eq:hatvphiNN0}
\lambda^{\mathstrut}_{N} p^{\mathstrut}_{N} \widehat{\varphi}^{N}_{N} \to 0
 \quad \textup{as  } N \to \infty.
\end{equation} 
One could use an asymptotic expansion to $\varphi^{\mathstrut}_{N}$, or just  set $\widehat{\varphi}^{N}_{N}$ to an arbitrary value, e.g.,  $0$. Note that $\varphi^{\mathstrut}_{N}$ does satisfy (\ref{eq:hatvphiNN0}), since, by Corollary \ref{cor:phiexp}(a), (\ref{eq:Tn}) and (\ref{eq:PCZ1is}),  
$\lambda^{\mathstrut}_{N} p^{\mathstrut}_{N} \varphi^{\mathstrut}_{N} = P_N (\zeta - Z_N) \to 0$.
Hence, condition (\ref{eq:hatvphiNN0}) is equivalent to 
\begin{equation}
\label{eq:hatphiNN}
\lambda^{\mathstrut}_{N} p^{\mathstrut}_{N} E_{\textup{abs}}(\widehat{\varphi}_{ N}^{ N}) \to 0 \quad \textup{as  } N \to \infty.
\end{equation}

Then, compute $\widehat{\varphi}^{N}_{N-1}$, \ldots, 
$\widehat{\varphi}^{N}_{0}$ by the backward recurrence  (cf.\ (\ref{eq:fnm1zbr}))

\begin{equation}
\label{eq:pephihatbr}
 \widehat{\varphi}^{N}_{n-1} = \frac{c_n - \widehat{\zeta}}{\mu_n} + \frac{\lambda_n}{\mu_n} \widehat{\varphi}_{n}^{ N}, \quad n = N, N-1, \ldots, 1.
\end{equation}
Note that, setting $z = \zeta$ in (\ref{eq:fnm1zbr}), the $\varphi_{n}$ satisfy 
\begin{equation}
\label{eq:epephihatbr}
\varphi_{n-1} = \frac{c_n - \zeta}{\mu_n} + \frac{\lambda_n}{\mu_n} \varphi_n, \quad n \in \mathbb{N}.
\end{equation}

As for the approximate relative costs 
$\widehat{b}_n^N$, they are computed by (cf.\ (\ref{eq:recbj}))
\begin{equation}
\label{eq:bhatnN}
\widehat{b}_0^N = b_0, \quad \widehat{b}_n^N = b_0 + \sum_{j=0}^{n-1} \widehat{\varphi}_j^N, \quad n = 1, \ldots, N.
\end{equation}

\subsection{Error analysis of backward recurrence scheme}
\label{s:beaacoi}
We next turn to analyzing the approximation errors of the computed $\widehat{\varphi}_{n}^{ N}$ and $\widehat{b}_n^N$.
We need the following preliminary result.

\begin{lemma}
\label{lma:fzbraps1}
For $n = 0, 1, \ldots, N,$
\begin{equation}
\label{eq:fnzbr}
\lambda_n p_n f_{n}(z) = \sum_{j=n+1}^ N c_j p_j -  z \sum_{j=n+1}^ N p_j + \lambda^{\mathstrut}_{N} p^{\mathstrut}_{N} f_N(z).
\end{equation}
\end{lemma}
\begin{proof}
We use backward induction on $n$. The case $n = N$ follows trivially. Suppose (\ref{eq:fnzbr}) holds for some $1 \leqslant n \leqslant N$. Then
\begin{align*}
\lambda_{n-1} p_{n-1} f_{n-1}(z) & = \lambda_{n-1} p_{n-1} \frac{c_n - z}{\mu_n} + \lambda_{n-1} p_{n-1} \frac{\lambda_n}{\mu_n} f_n(z) \\
& = p_{n} (c_n - z) +   \lambda_n p_{n} f_n(z) \\
& = p_{n} (c_n - z) +   \sum_{j=n+1}^ N c_j p_j -  z \sum_{j=n+1}^ N p_j + \lambda^{\mathstrut}_{N} p^{\mathstrut}_{N} f_N(z) \\
& = \sum_{j=n}^ N c_j p_j -  z \sum_{j=n}^ N p_j + \lambda^{\mathstrut}_{N} p^{\mathstrut}_{N} f_N(z),
\end{align*}
using in turn (\ref{eq:fnm1zbr}), 
$\lambda_{n-1} p_{n-1} = \mu_n p_n$ and (\ref{eq:fnzbr}), 
which completes the induction.
\end{proof}

The next result follows readily from Lemma \ref{lma:fzbraps1} by taking $z = \zeta$ and $z = \widehat{\zeta}$.

\begin{corollary}
\label{lma:braps1}
For $n = 0, 1, \ldots, N-1,$
\begin{itemize}
\item[\textup{(a)}] $\displaystyle 
\lambda_n p_n \varphi_{n} = \sum_{j=n+1}^ N c_j p_j -  \zeta \sum_{j=n+1}^ N p_j + \lambda^{\mathstrut}_{N} p^{\mathstrut}_{N} \varphi_N;$
\item[\textup{(b)}] $\displaystyle \lambda_n p_n \widehat{\varphi}_{n}^N = \sum_{j=n+1}^ N c_j p_j -  \widehat{\zeta} \sum_{j=n+1}^ N p_j + \lambda^{\mathstrut}_{N} p^{\mathstrut}_{N} \widehat{\varphi}_{ N}^{ N}.$
\end{itemize}
\end{corollary}

The next result gives expressions for  the errors
$E_{\textup{abs}}(\widehat{\varphi}_{n}^N)$ and $E_{\textup{abs}}(\widehat{b}_n^N)$.
 
\begin{lemma}
\label{lma:bcphisdif}
\textup{ }
\begin{itemize}
\item[\textup{(a)}] $\displaystyle
\lambda_n p_n E_{\textup{abs}}(\widehat{\varphi}_{n}^N) =   \lambda^{\mathstrut}_{N} p^{\mathstrut}_{N} E_{\textup{abs}}(\widehat{\varphi}_{ N}^{ N}) -   (P_{N} - P_n) E_{\textup{abs}}(\widehat{\zeta}),$ \, $0 \leqslant n < N;$
\item[\textup{(b)}] $\displaystyle E_{\textup{abs}}(\widehat{b}_n^N) =   \lambda^{\mathstrut}_{N} p^{\mathstrut}_{N} E_{\textup{abs}}(\widehat{\varphi}_{ N}^{ N}) \sum_{j=0}^{n-1}  \frac{1}{\lambda_j p_j} -   E_{\textup{abs}}(\widehat{\zeta})  \sum_{j=0}^{n-1} \frac{P_{N} - P_j}{\lambda_j p_j},$ $1 \leqslant  n \leqslant N.$
\end{itemize}
\end{lemma}
\begin{proof}
(a) This part follows straightforwardly from Corollary \ref{lma:braps1}, by subtracting the expressions given there for $\lambda_n p_n \widehat{\varphi}_{n}^N$ and $\lambda_n p_n \varphi_{n}$.

(b) Using (\ref{eq:recbj}), (\ref{eq:bhatnN}) and part (a), we obtain
\begin{align*}
E_{\textup{abs}}(\widehat{b}_n^N) & = \sum_{j=0}^{n-1} E_{\textup{abs}}(\widehat{\varphi}_j^N)  = \lambda^{\mathstrut}_{N} p^{\mathstrut}_{N} E_{\textup{abs}}(\widehat{\varphi}_{ N}^{ N}) \sum_{j=0}^{n-1}  \frac{1}{\lambda_j p_j} -   E_{\textup{abs}}(\widehat{\zeta})  \sum_{j=0}^{n-1} \frac{P_N - P_j}{\lambda_j p_j}.
\end{align*}
\end{proof}

The following result is the counterpart of 
Proposition \ref{pro:phihatsol} for the backward recurrence scheme.
It relates the approximation errors 
in the computed
$\widehat{\varphi}_{n}^N$ and $\widehat{b}_n^N$ to those  in the input $\widehat{\zeta}$, by identifying the corresponding \emph{asymptotic error amplificaton factors} as $N \to \infty$.
 
\begin{proposition}[Asymptotic error amplification]
\label{pro:brphihatsol}
For each $n,$ as $N \to \infty,$
\begin{itemize}
\item[\textup{(a)}]
$\displaystyle E_{\textup{abs}}(\widehat{\varphi}_{n}^N) \to -T_{n+1}^- E_{\textup{abs}}(\widehat{\zeta});$
\item[\textup{(b)}]
$\displaystyle E_{\textup{rel}}(\widehat{\varphi}_{n}^N) \to   -   \zeta \frac{T_{n+1}^-}{\varphi_n} E_{\textup{rel}}(\widehat{\zeta}) = \frac{\zeta }{Z_{n+1}^- - \zeta} E_{\textup{rel}}(\widehat{\zeta});$ 
\item[\textup{(c)}] $\displaystyle E_{\textup{abs}}(\widehat{b}_n^N) \to  -    T_{n0} E_{\textup{abs}}(\widehat{\zeta}) ;$
\item[\textup{(d)}]
$\displaystyle E_{\textup{rel}}(\widehat{b}_{n}^N) \to
- \zeta  \frac{T_{n0}}{b_n} \, E_{\textup{rel}}(\widehat{\zeta}) = 
-   \frac{\zeta}{{b_0}/{T_{n0}} + Z_{n0} - \zeta} \, E_{\textup{rel}}(\widehat{\zeta}).$ 
\end{itemize}
\end{proposition}
\begin{proof}
(a) From Lemma \ref{lma:bcphisdif}(a), (\ref{eq:hatvphiNN0}), (\ref{eq:hatphiNN}), Lemma \ref{lma:CPm}(a) and (\ref{eq:Tnmkeil}), we obtain
\begin{align*}
E_{\textup{abs}}(\widehat{\varphi}_{n}^N) & = \frac{\lambda^{\mathstrut}_{N} p^{\mathstrut}_{N}}{\lambda_n p_n} E_{\textup{abs}}(\widehat{\varphi}_{ N}^{ N}) -  \frac{P_N - P_n}{\lambda_n p_n} E_{\textup{abs}}(\widehat{\zeta}) \\& \to -  \frac{\bar{P}_n}{\lambda_n p_n} E_{\textup{abs}}(\widehat{\zeta}) = -T_{n+1}^- E_{\textup{abs}}(\widehat{\zeta}) \textup{  as  } N \to \infty.
\end{align*}

(b) This part follows from part (a) and Proposition \ref{pro:bfzexp}(a).

(c) This part follows by letting $N \to \infty$ in Lemma \ref{lma:bcphisdif}(b) and using (\ref{eq:Tn0}).

(d) The result follows using part (c) and Proposition \ref{pro:bfzexp}(b).
\end{proof}

Thus, Proposition \ref{pro:brphihatsol}(a, c) ensures that, for a fixed state $n$, the computed $\widehat{\varphi}_{n}^N$ and $\widehat{b}_{n}^N$ are asymptotically related to the exact $\varphi_n$ and  $b_n$ by
\begin{equation}
\label{eq:brasympphib}
\widehat{\varphi}_{n}^N \to \varphi_n -T_{n+1}^- E_{\textup{abs}}(\widehat{\zeta})
\quad \textup{and} \quad 
\widehat{b}_{n}^N \to b_n - T_{n0} E_{\textup{abs}}(\widehat{\zeta}) \enspace
\textup{as} \enspace N \to \infty.
\end{equation}

How does this compare to the quality of approximations $\widehat{\varphi}_{n}$ and $\widehat{b}_{n}$ computed by forward recursion in \S\ref{s:fraeas}?
The next result shows that the backward recurrence approximations $\widehat{\varphi}_{n}^N$ and $\widehat{b}_{n}^N$ are asymptotically much more accurate, in that $E_{\textup{abs}}(\widehat{\varphi}_{n}^N)/E_{\textup{abs}}(\widehat{\varphi}_{n}) \approx 0$ and $E_{\textup{abs}}(\widehat{b}_{n}^N)/E_{\textup{abs}}(\widehat{b}_{n}) \approx 0$ for large $n$ and $N > n$.

\begin{proposition}[Relative accuracy of backward versus forward recurrence]
\label{pro:compfrbr} \textup{ }
\begin{itemize}
\item[\textup{(a)}] $\displaystyle \lim_{n \to \infty} \lim_{N \to \infty} E_{\textup{abs}}(\widehat{\varphi}_{n}^N)/E_{\textup{abs}}(\widehat{\varphi}_{n}) = 0;$
\item[\textup{(b)}] $\displaystyle \lim_{n \to \infty} \lim_{N \to \infty} E_{\textup{abs}}(\widehat{b}_{n}^N)/E_{\textup{abs}}(\widehat{b}_{n}) = 0.$
\end{itemize}
\end{proposition}
\begin{proof}
(a) From Propositions \ref{pro:phihatsol}(a) and \ref{pro:brphihatsol}(a), and the leftmost identity in (\ref{eq:Tnmkeil}) we have, for each $n$,
\[
\lim_{N \to \infty} \frac{E_{\textup{abs}}(\widehat{\varphi}_{n}^N)}{E_{\textup{abs}}(\widehat{\varphi}_{n})} = -\frac{T_{n+1}^-}{T_n^+} = -\frac{\bar{P}_{n}}{P_n}.
\]
The result now follows since ${\bar{P}_{n}}/{P_n} \searrow 0$ as $n \to \infty$.

(b) From Propositions \ref{pro:phihatsol}(c) and \ref{pro:brphihatsol}(c),  we have, for each $n$,
\[
\lim_{N \to \infty} \frac{E_{\textup{abs}}(\widehat{b}_{n}^N)}{E_{\textup{abs}}(\widehat{b}_{n})} = -\frac{T_{n0}}{T_{0n}}.
\]
The result now follows since ${T_{n0}}/{T_{0n}} \searrow 0$ as $n \to \infty$ by  Lemma \ref{lma:ratioTn00n}(b).
\end{proof}

\section{A mixed forward--backward recurrence scheme}
\label{s:afbr}
The results in \S\ref{s:beaacoi} show that 
backward recurrence gives more accurate computed solutions 
to Poisson's equation than forward recurrence for 
 large states, whereas it gives less accurate results for small states.
This insight leads us to propose a novel mixed \emph{forward--backward recurrence} scheme that keeps the pros and avoids the cons of the pure schemes.

In light of Propositions \ref{pro:phihatsol} and \ref{pro:brphihatsol},  let us define, for a given large integer $N$ and for states $n < N$, 
 the approximate marginal relative costs 
\[
\widetilde{\varphi}_{n}^{ N} \triangleq 
\begin{cases}
\widehat{\varphi}_n, & \textup{ if }  P_n \leqslant  \bar{P}_n \\
\widehat{\varphi}_{n}^{ N}, & \textup{ otherwise,}
\end{cases}
\]
noting that $P_n \leqslant  \bar{P}_n$ if and only if $P_n \leqslant 1/2$,
and the approximate relative costs
\[
\widetilde{b}_n^N \triangleq 
\begin{cases}
\widehat{b}_n, & \textup{ if }  T_{0n} \leqslant  T_{n0} \\
\widehat{b}_n^N, & \textup{ otherwise.}
\end{cases}
\]

\begin{remark}
\label{re:bothsdc0}
Since both $\bar{P}_n/P_n$ and $T_{n0}/T_{0n}$ strictly decrease to $0$ as $n \to \infty$  (see Lemma \ref{lma:ratioTn00n}),  letting $m \triangleq \min \{n~\geqslant~0\colon  \bar{P}_n/P_n < 1\}$ and $M \triangleq \min \{n \geqslant 1\colon T_{n0}/T_{0n} < 1\}$, we have
$P_n \leqslant  \bar{P}_n$ if and only if $n < m$, and $T_{0n}  \leqslant T_{n0}$ if and only if $n < M$. Note that it must be the case that $m \leqslant M$ by Lemma \ref{lma:ratioTn00n}(a).
\end{remark}

Now, letting $m$ and $M$ be as in Remark \ref{re:bothsdc0}, define 
\[
A_n \triangleq 
\begin{cases}
T_n^+, & \textup{ if } n <  m \\
-T_{n+1}^-, & \textup{ otherwise,}
\end{cases}
\quad
\textup{and}
\quad
B_n \triangleq 
\begin{cases}
T_{0n}, & \textup{ if } n <  M \\
-T_{n0}, & \textup{ otherwise.}
\end{cases}
\]

In the following result,  we assume that relative errors are well defined.

\begin{proposition}[Asymptotic error amplification]
\label{pro:fbrphihatsol}
For each $n,$ as $N \to \infty,$
\begin{itemize}
\item[\textup{(a)}]
$\displaystyle E_{\textup{abs}}(\widetilde{\varphi}_{n}^N) \to A_{n} E_{\textup{abs}}(\widehat{\zeta});$
\item[\textup{(b)}]
$\displaystyle E_{\textup{rel}}(\widetilde{\varphi}_{n}^N) \to   \frac{\zeta A_{n}}{\varphi_n} E_{\textup{rel}}(\widehat{\zeta});$ 
\item[\textup{(c)}] $\displaystyle E_{\textup{abs}}(\widetilde{b}_n^N) \to  B_n E_{\textup{abs}}(\widehat{\zeta});$
\item[\textup{(d)}]
$\displaystyle E_{\textup{rel}}(\widetilde{b}_{n}^N) \to
\frac{\zeta B_n}{b_n} \, E_{\textup{rel}}(\widehat{\zeta}).$
\end{itemize}
\end{proposition}
\begin{proof}
The result follows directly from Propositions \ref{pro:phihatsol} and  \ref{pro:brphihatsol}.
\end{proof}

The following result shows that forward--backward recurrence is asymptotically more accurate than either of the pure recurrence schemes.
\begin{corollary}
\label{cor:fbrphihatsol} For each $n,$
\begin{itemize}
\item[\textup{(a)}]
$\displaystyle \lim_{N \to \infty} |E_{\textup{abs}}(\widetilde{\varphi}_{n}^N)| = 
\min\{|E_{\textup{abs}}(\widehat{\varphi}_{n})|, \lim_{N \to \infty} |E_{\textup{abs}}(\widehat{\varphi}_{n}^N)|\};$
\item[\textup{(b)}]
 $\displaystyle \lim_{N \to \infty} |E_{\textup{abs}}(\widetilde{b}_n^N)| =  \min\{|E_{\textup{abs}}(\widehat{b}_{n})|, \lim_{N \to \infty} |E_{\textup{abs}}(\widehat{b}_{n}^N)|\}.$
\end{itemize}
\end{corollary}
\begin{proof}
The result follows directly from Propositions \ref{pro:phihatsol},  \ref{pro:brphihatsol} and \ref{pro:fbrphihatsol}.
\end{proof}

\subsection{Error analysis of bias and asymptotic variance computation}
\label{s:fbeacrvav}
Let
$\beta_{0, N}$ and 
$\sigma_N^2$ be as in \S\ref{s:eacrvav}, which we assume satisfy $\beta_{0, N} \to \beta_0 = 
- \sum_{n=0}^\infty \bar{P}_n \varphi_n$ and $\sigma_N^2 \to \sigma^2 = 2 \sum_{n=0}^\infty \lambda_n p_n \varphi_n^2$ as $N \to \infty$ (cf. Proposition \ref{pro:betasigma2phi}), 
and consider the computed approximations
$\widetilde{\beta}_{0, N} \triangleq 
-\sum_{n =0}^{N} \bar{P}_n \widetilde{\varphi}_{n}^N$ and
$\widetilde{\sigma}_N^2 \triangleq  2 \sum_{n =0}^{N} \lambda_n p_n (\widetilde{\varphi}_{n}^{ N})^2$.
The next result shows that, unlike the $\widehat{\beta}_{0, N}$ and $\widehat{\sigma}_N^2$ in \S\ref{s:eacrvav}, $\widetilde{\beta}_{0, N}$ and $\widetilde{\sigma}_N^2$ have bounded approximation errors as $N \to \infty$, for which explicit expressions are given, 
provided that $T_{p 0} < \infty$ (see Remark \ref{re:T0np1div}(c)) and   $\widehat{\varphi}_{ N}^{ N}$ 
satisfies more stringent asymptotic accuracy requirements than (\ref{eq:hatphiNN}). 

We consider the following conditions on the quality of  approximation $\widehat{\varphi}_{ N}^{ N}$: 
\begin{equation}
\label{eq:asbetatild}
\lambda^{\mathstrut}_{N} p^{\mathstrut}_{N}  T_{N0} \, E_{\textup{abs}}(\widehat{\varphi}_{ N}^{ N}) \to 0 \enspace \textup{as} \enspace N \to \infty,
\end{equation}
\begin{equation}
\label{eq:Hasbetatild}
\lambda^{\mathstrut}_{N} p^{\mathstrut}_{N}  H_{N0} \, E_{\textup{abs}}(\widehat{\varphi}_{ N}^{ N}) \to 0 \enspace \textup{as} \enspace N \to \infty,
\end{equation}
and
\begin{equation}
\label{eq:asbetatild2}
\lambda^{\mathstrut}_{N} p^{\mathstrut}_{N}  \sqrt{T_{0N}} \, E_{\textup{abs}}(\widehat{\varphi}_{ N}^{ N}) \to 0 \enspace \textup{as} \enspace N \to \infty.
\end{equation}

\begin{proposition}
\label{pro:fbnbetasigma2phi} Suppose that $T_{p0} < \infty.$ Then$,$
\begin{itemize}
\item[\textup{(a)}] under $(\ref{eq:asbetatild}),$ $\displaystyle E_{\textup{abs}}(\widetilde{\beta}_{0, N}) \to   (T_{p0}-T_{m 0}) E_{\textup{abs}}(\widehat{\zeta})$ as $N \to \infty;$ 
\item[\textup{(b)}] under $(\ref{eq:asbetatild})$, $(\ref{eq:Hasbetatild})$ and $(\ref{eq:asbetatild2}),$
\[
E_{\textup{abs}}(\widetilde{\sigma}_N^2) \to   
 4 \beta_{m} E_{\textup{abs}}(\widehat{\zeta})  + 
2 (T_{p0} + T_{0 m} - T_{m0}) E_{\textup{abs}}^2(\widehat{\zeta})
\enspace \textup{as} \enspace N \to \infty. 
\]
\end{itemize}
\end{proposition}
\begin{proof}
(a) From Proposition \ref{pro:phihatsol}(a) and Lemma \ref{lma:bcphisdif}(a) we have, for $N \geqslant m+2$, 
\begin{align*}
E_{\textup{abs}}(\widetilde{\beta}_{0, N}) & = - \sum_{n=0}^{N-1} \bar{P}_n E_{\textup{abs}}(\widetilde{\varphi}_{n}^{ N}) =
-\sum_{n=0}^{m-1} \bar{P}_n E_{\textup{abs}}(\widehat{\varphi}_n ) - \sum_{n=m}^{N-1} \bar{P}_n E_{\textup{abs}}(\widehat{\varphi}_{n}^{ N})
\\
& = -E_{\textup{abs}}(\widehat{\zeta}) \sum_{n=0}^{m-1} \bar{P}_n T_n^+  - \lambda^{\mathstrut}_{N} p^{\mathstrut}_{N} E_{\textup{abs}}(\widehat{\varphi}_{ N}^{ N}) \sum_{n=m}^{N-1} \frac{\bar{P}_n}{\lambda_n p_n}  \\
& \qquad +  E_{\textup{abs}}(\widehat{\zeta}) \sum_{n=m}^{N-1} \frac{\bar{P}_n}{\lambda_n p_n} (P_N - P_n).
\end{align*}

Now, from the above identity we obtain, as $N \to \infty$,
 \begin{align*}
E_{\textup{abs}}(\widetilde{\beta}_{0, N})
& \to 
 \bigg(\sum_{n=m}^\infty \frac{\bar{P}_n^2}{\lambda_n p_n}-\sum_{n=0}^{m-1}  \frac{P_n \bar{P}_n}{\lambda_n p_n}\bigg) E_{\textup{abs}}(\widehat{\zeta}) \\
& = \bigg(\sum_{n=0}^\infty \frac{\bar{P}_n^2}{\lambda_n p_n}-\sum_{n=0}^{m-1}  \frac{\bar{P}_n}{\lambda_n p_n}\bigg) E_{\textup{abs}}(\widehat{\zeta}) = (T_{p0}-T_{m 0}) E_{\textup{abs}}(\widehat{\zeta}),
\end{align*}
where we have used (\ref{eq:Tn}), (\ref{eq:Tn0}), (\ref{eq:asbetatild}),  and (\ref{eq:Tp0series}).

(b)
We have 
\begin{align*}
\widetilde{\sigma}_{N}^2 & = 
2 \sum_{n=0}^{N-1} \lambda_n p_n 
\big(\widetilde{\varphi}_{ n}^{ N}\big)^2 =
2 \sum_{n=0}^{m-1} \lambda_n p_n 
\widehat{\varphi}_n^2 +
2 \sum_{n=m}^{N-1} \lambda_n p_n 
\big(\widehat{\varphi}_{ n}^{ N}\big)^2.
\end{align*}

Now, on the one hand, using Proposition \ref{pro:phihatsol}(a) and (\ref{eq:Tn}) gives
\begin{align*}
\sum_{n=0}^{m-1} \lambda_n p_n 
\widehat{\varphi}_n^2 & = 
\sum_{n=0}^{m-1} \lambda_n p_n 
\big(\varphi_n + T_n^+ E_{\textup{abs}}(\widehat{\zeta})\big)^2 \\
& = \sum_{n=0}^{m-1} \lambda_n p_n \varphi_n^2 
+ 2 E_{\textup{abs}}(\widehat{\zeta}) \sum_{n=0}^{m-1} P_n \varphi_n + 
E_{\textup{abs}}^2(\widehat{\zeta}) \sum_{n=0}^{m-1} \frac{P_n^2}{\lambda_n p_n}.
\end{align*}

On the other hand, using Lemma \ref{lma:bcphisdif}(a),  $\sum_{n=m}^{N-1} \lambda_n p_n 
\big(\widehat{\varphi}_{ n}^{ N}\big)^2$ equals
\begin{align*}
 &
\sum_{n=m}^{N-1} \lambda_n p_n 
\Big(\varphi_n + \frac{\lambda^{\mathstrut}_N p^{\mathstrut}_N}{\lambda_n p_n} E_{\textup{abs}}(\widehat{\varphi}_{ N}^{ N}) -   \frac{P_N - P_n}{\lambda_n p_n} E_{\textup{abs}}(\widehat{\zeta})\Big)^2 \\
& = \sum_{n=m}^{N-1} \lambda_n p_n \varphi_n^2 + 
2 \sum_{n=m}^{N-1} \varphi_n \big(\lambda^{\mathstrut}_N p^{\mathstrut}_N E_{\textup{abs}}(\widehat{\varphi}_{ N}^{ N}) -   E_{\textup{abs}}(\widehat{\zeta}) (P_N - P_n)\big) \\
& \quad + \sum_{n=m}^{N-1} \frac{1}{\lambda_n p_n} \big(\lambda^{\mathstrut}_N p^{\mathstrut}_N E_{\textup{abs}}(\widehat{\varphi}_{ N}^{ N}) -   (P_N - P_n) E_{\textup{abs}}(\widehat{\zeta})\big)^2 \\
& = \sum_{n=m}^{N-1} \lambda_n p_n \varphi_n^2 + 
2 \lambda^{\mathstrut}_N p^{\mathstrut}_N E_{\textup{abs}}(\widehat{\varphi}_{ N}^{ N}) \sum_{n=m}^{N-1} \varphi_n   \\
& \quad -   2 E_{\textup{abs}}(\widehat{\zeta}) \sum_{n=m}^{N-1}(P_N - P_n) \varphi_n + \big(\lambda^{\mathstrut}_N p^{\mathstrut}_N E_{\textup{abs}}(\widehat{\varphi}_{ N}^{ N})\big)^2 \sum_{n=m}^{N-1} \frac{1}{\lambda_n p_n} \\
& \quad
 -  2 E_{\textup{abs}}(\widehat{\zeta}) \lambda^{\mathstrut}_N p^{\mathstrut}_N E_{\textup{abs}}(\widehat{\varphi}_{ N}^{ N}) \sum_{n=m}^{N-1} \frac{P_N - P_n}{\lambda_n p_n}   + E_{\textup{abs}}^2(\widehat{\zeta}) \sum_{n=m}^{N-1} \frac{(P_N - P_n)^2}{\lambda_n p_n},
\end{align*}
and hence we have, as $N \to \infty$,  
\[
\sum_{n=m}^{N-1} \lambda_n p_n 
\big(\widehat{\varphi}_{ n}^{ N}\big)^2 \to \sum_{n=m}^{\infty} \lambda_n p_n \varphi_n^2 -   2 E_{\textup{abs}}(\widehat{\zeta}) \sum_{n=m}^{\infty} \bar{P}_n \varphi_n
 + E_{\textup{abs}}^2(\widehat{\zeta}) \sum_{n=m}^{\infty} \frac{\bar{P}_n^2}{\lambda_n p_n},
\]
as the other terms vanish under the assumptions. 
Thus, 
\begin{align*}
\lambda^{\mathstrut}_N p^{\mathstrut}_N E_{\textup{abs}}(\widehat{\varphi}_{ N}^{ N}) \sum_{n=m}^{N-1} \varphi_n & = 
\lambda^{\mathstrut}_N p^{\mathstrut}_N (\beta_N - \beta_m) E_{\textup{abs}}(\widehat{\varphi}_{ N}^{ N}) \\
& \approx
\lambda^{\mathstrut}_N p^{\mathstrut}_N (H_{N0} - \zeta  T_{N0}) E_{\textup{abs}}(\widehat{\varphi}_{ N}^{ N}) \to 0 \enspace \textup{as} \enspace N \to \infty,
\end{align*}
where we have used Propositions \ref{pro:betasigma2phi}(a) and
\ref{pro:bfzexp}(b), and then (\ref{eq:asbetatild}) and (\ref{eq:Hasbetatild}).

Also, as $N \to \infty$, 
\begin{align*}
\big(\lambda^{\mathstrut}_N p^{\mathstrut}_N E_{\textup{abs}}(\widehat{\varphi}_{ N}^{ N})\big)^2 \sum_{n=m}^{N-1} \frac{1}{\lambda_n p_n} & \approx \big(\lambda^{\mathstrut}_N p^{\mathstrut}_N E_{\textup{abs}}(\widehat{\varphi}_{ N}^{ N})\big)^2 (T_{0N} + T_{N0})  \\
& = \Big(\lambda^{\mathstrut}_N p^{\mathstrut}_N E_{\textup{abs}}(\widehat{\varphi}_{ N}^{ N}) \sqrt{T_{0N} + T_{N0}}\Big)^2
 \to 0,
\end{align*}
where we have used (\ref{eq:T0np1H0np1}), (\ref{eq:Tn0}), (\ref{eq:asbetatild}) and  (\ref{eq:asbetatild2}). 

Further, from (\ref{eq:Tn0}) and (\ref{eq:asbetatild}) we have, as $N \to \infty$, 
\begin{align*}
\lambda^{\mathstrut}_N p^{\mathstrut}_N E_{\textup{abs}}(\widehat{\varphi}_{ N}^{ N}) \sum_{n=m}^{N-1} \frac{P_N - P_n}{\lambda_n p_n} & \approx \lambda^{\mathstrut}_N p^{\mathstrut}_N T_{N0} E_{\textup{abs}}(\widehat{\varphi}_{ N}^{ N})  \to 0.
\end{align*}

It follows from the above that, as $N \to \infty$,
$E_{\textup{abs}}(\widetilde{\sigma}_N^2)$ converges to
\begin{align*}
& 4  \Bigg(\sum_{n=0}^{m-1} P_n \varphi_n - \sum_{n=m}^{\infty} \bar{P}_n \varphi_n\Bigg) E_{\textup{abs}}(\widehat{\zeta}) + 
 2 \Bigg(\sum_{n=0}^{m-1} \frac{P_n^2}{\lambda_n p_n}  
 + \sum_{n=m}^{\infty} \frac{\bar{P}_n^2}{\lambda_n p_n}\Bigg) E_{\textup{abs}}^2(\widehat{\zeta}) \\
& = 4  \Bigg(\sum_{n=0}^{m-1} \varphi_n - \sum_{n=0}^{\infty} \bar{P}_n \varphi_n\Bigg) E_{\textup{abs}}(\widehat{\zeta}) + 
 2 \Bigg(\sum_{n=0}^{m-1} \frac{P_n^2-\bar{P}_n^2}{\lambda_n p_n}  
 + \sum_{n=0}^{\infty} \frac{\bar{P}_n^2}{\lambda_n p_n}\Bigg) E_{\textup{abs}}^2(\widehat{\zeta}) \\
& = 4  \beta_m E_{\textup{abs}}(\widehat{\zeta}) + 
 2 (T_{0m} - T_{m0}   + T_{p0}) E_{\textup{abs}}^2(\widehat{\zeta}),
\end{align*}
using Proposition \ref{pro:betasigma2phi}(a), (\ref{eq:T0np1H0np1}), (\ref{eq:Tn0}), $P_n^2-\bar{P}_n^2 = P_n-\bar{P}_n$, and (\ref{eq:Tp0series}).
\end{proof}

\section{An example}
\label{s:anex}
This section illustrates the application of the above results to the 
 M/M/1+M queueing model with deadlines to the end of service considered in \S \ref{s:intro}.

For this model the steady-state probabilities are given by (see \cite[p.\ 89]{anckGaf62})
\begin{equation}
\label{eq:piimm1m}
p_n = \frac{e^{-\kappa} }{\mathcal{P}(\alpha, \kappa)} \frac{\kappa^{ \alpha+ n}}{\Gamma( \alpha+ n +1)}, \quad 
n \in \mathbb{N}_0,
\end{equation}
where $\alpha \triangleq \mu/\theta$ and $\kappa \triangleq \lambda/\theta$.
Note that $\Gamma(a)$ is the gamma function, $\gamma(a, x) \triangleq \int_0^x t^{a-1} e^{-t} \, dt$ and $\Gamma(a, x) \triangleq \int_x^\infty t^{a-1} e^{-t} \, dt$ are the  lower and upper incomplete gamma functions, and $\mathcal{P}(a, x) \triangleq \gamma(a, x)/\Gamma(a)$ and  
$\mathcal{Q}(a, x) \triangleq \Gamma(a, x)/\Gamma(a)$ are the lower and  upper normalized gamma functions, respectively.  See \cite[\S8.2]{NIST:DLMF}.

The mean steady-state cost has the evaluation
\begin{equation}
\label{eq:zetamm1mpm}
\zeta = \lambda - \mu (1-p_0) = 
 \lambda - \mu \bigg(1- 
    \frac{\theta \kappa^{\alpha} e^{-\kappa}}{\mu \gamma(\alpha,\kappa)}\bigg) = 
\lambda - \mu + 
    \frac{\theta \kappa^{\alpha} e^{-\kappa}}{\gamma(\alpha,\kappa)},    
\end{equation}
and the cumulative steady-state probabilities 
are given by
\begin{equation}
\label{eq:MM1MPi}
P_n = 1 - \frac{\mathcal{P}( \alpha+ n +1,\kappa)}{\mathcal{P}(\alpha ,\kappa)}, \quad 
n \in \mathbb{N}_0.
\end{equation}

From (\ref{eq:MM1MPi}), expressions for 
the mean first-passage times and costs considered in the above analyses are readily obtained.
From these, and using Corollary \ref{cor:phiexp}, the following analytical expression for the marginal relative cost is obtained:
\begin{equation}
\label{eq:phinexact}
\varphi_n =  1 - 
    \frac{\gamma( \alpha+ n +1,\kappa)}{\gamma(\alpha,\kappa) \, \kappa
   ^{n + 1}}  = 1 - \frac{\kappa^{\alpha-1}}{\gamma(\alpha,\kappa)} \int_0^\kappa (t/\kappa)^{ \alpha+ n} e^{-t} \, dt, \quad 
n \in \mathbb{N}_0.
\end{equation}
Now, it is evident from the rightmost expression in (\ref{eq:phinexact}) that $\varphi_n \nearrow 1$ as $n \to \infty$.

Recall that Table \ref{tab:numex1} in \S\ref{s:intro} illustrated 
the explosive numerical instability of the approximate solution to Poisson's equation by standard forward recurrence. 
We now consider the same instance, but approximate 
the marginal costs $\varphi_n$ by the $\widetilde{\varphi}_n^N$ calculated by the forward--backward recurrence scheme presented in \S\ref{s:abr}.
To test the accuracy of results, the values obtained were compared to the computed values of the $\varphi_n$ using (\ref{eq:phinexact}), which we denote by $fl(\varphi_n)$, where $fl(x)$ denotes the floating-point approximation of a number $x$. 
Computations were done in Matlab with standard double-precision arithmetic, so the relative error of $fl(x)$ is bounded above by the 
unit roundoff $u = 2^{-53} \approx 1.1 \times 10^{-16}$.
We thus consider that the approximation to $\zeta$ computed by Matlab is 
$\widehat{\zeta} = fl(\zeta)$.

We found that taking $N = 42$ and $\widehat{\varphi}_N^N = 0$ suffices to obtain extremely accurate approximations to $\varphi_n$ for the values of $n$ in Table  \ref{tab:numex1} in which inaccuracies were evident, i.e., $n$ larger than $11$. Table \ref{tab:numex2} shows the results. 
The $\widetilde{\varphi}_{n}^N$ column shows the evaluations of such quantities with 15 significant digits, which precisely match those of 
the $fl(\varphi_n)$.
The column labeled $\zeta A_{n}/\varphi_n$ evaluates the theoretical
 asymptotic relative-error amplification factors in Proposition \ref{pro:fbrphihatsol}(b). 
The results there show that the relative error of $\widetilde{\varphi}_{n}^N$ is, in theory, substantially reduced with respect to that of $\widehat{\zeta}$. 
The next column, labeled $2^{53} |E_{\textup{rel}}(\widetilde{\varphi}_{n}^N)|$, evaluates the ratio of the absolute value of the actual relative error of $\widetilde{\varphi}_{n}^N$ (evaluated as $|\widetilde{\varphi}_{n}^N - fl(\varphi_n)|/fl(\varphi_n)$) to that of $\widehat{\zeta}$, which is taken equal to  the unit roundoff $u$.
The results differ slightly from the theory, but they show that the relative error of $\widetilde{\varphi}_{n}^N$ is near $u$ in the worst case.

The values shown in the last two columns, labeled $T_{n+1}^-$ and $T_{n}^+$, explain the vastly improved accuracy of forward--backward recurrence  with respect to forward recurrence.
Recall from Proposition \ref{pro:fbrphihatsol}(a) that the absolute approximation error of $\widetilde{\varphi}_{n}^N$, for large $n$ and $N$, is approximately proportional to $T_{n+1}^-$, while Proposition \ref{pro:phihatsol}(a) shows that the absolute approximation error of $\widehat{\varphi}_{n}$ is proportional to $T_{n}^+$.
These columns show that the $T_{n+1}^-$ are very small and vanish as $n$ grows, while the $T_{n}^+$ quickly grow to infinity.

\begin{table}[tbh!]
\caption{Accurate numerical computation of $\varphi_n$ by forward--backward recurrence.} \label{tab:numex2}
\scriptsize 
 \vspace{.1in}
\begin{center}
\begin{tabular}{|cccccc|}
$n$ & $\widetilde{\varphi}_{n}^N$ & $\zeta A_{n}/\varphi_n$ & \multicolumn{1}{c}{$2^{53} |E_{\textup{rel}}(\widetilde{\varphi}_{n}^N)|$} & $T_{n+1}^-$ & $T_{n}^+$ \\ \hline
$12$  & $0.925174342237504$ & $6.47 \times 10^{-2}$ & $0$ & $0.150$ & $8.4 \times 10^7$ \\
$13$  & $0.930359089413224$ & $6.00 \times 10^{-2}$ & $0$ & $0.140$ & $7.0 \times 10^8$ \\
$14$  & $0.934875921107126$ & $5.57 \times 10^{-2}$ & $1.07$ & $0.131$ & $6.2 \times 10^9$ \\
$15$  & $0.938845492334662$ & $5.21 \times 10^{-2}$  & $0$ & $0.123$ & $5.9 \times 10^{10}$  \\
$16$  & $0.942361160780650$ & $4.89 \times 10^{-2}$ & $0$ & $0.116$ & $5.9 \times 10^{11}$ \\
$17$  & $0.945496267896444$ & $4.61 \times 10^{-2}$ & $1.06$ & $0.109$ & $6.2 \times 10^{12}$ \\
$18$  & $0.948309214061184$ & $4.36 \times 10^{-2}$ & $1.05$ & $0.104$ & $6.9 \times 10^{13}$ \\
$19$  & $0.950847068147842$ & $4.13 \times 10^{-2}$ & $1.05$ & $0.099$ & $8.4 \times 10^{14}$ \\
$20$  & $0.953148181463212$ & $3.93 \times 10^{-2}$ & $0$ & $0.094$ & $9.8 \times 10^{15}$ \\
$21$  & $0.955244111686174$ & $3.75 \times 10^{-2}$ & $0$ & $0.090$ & $1.3 \times 10^{17}$ \\
$22$  & $0.957161059916347$ & $3.58 \times 10^{-2}$ & $1.04$ & $0.086$ & $1.7 \times 10^{18}$ \\
$23$  & $0.958920958494403$ & $3.42 \times 10^{-2}$ & $1.04$ & $0.082$ & $2.3 \times 10^{19}$ \\
$24$  & $0.960542304575400$ & $3.28 \times 10^{-2}$ & $0$ & $0.079$ & $3.4 \times 10^{20}$ \\
$25$  & $0.962040806065022$ & $3.15 \times 10^{-2}$ & $1.04$ & $0.076$ & $5.0 \times 10^{21}$ \\
$26$  & $0.963429887334373$ & $3.03 \times 10^{-2}$ & $0$ & $0.073$ & $7.8 \times 10^{22}$ \\
$27$  & $0.964721088932251$ & $2.92 \times 10^{-2}$ & $1.04$ & $0.071$ & $1.3 \times 10^{24}$ \\
$28$  & $0.965924386304869$ & $2.82 \times 10^{-2}$ & $0$ & $0.068$ & $2.1 \times 10^{25}$ \\
$29$  & $0.967048446017873$ & $2.72 \times 10^{-2}$ & $0$ & $0.066$ & $3.6 \times 10^{26}$ \\ \hline
\end{tabular}
\end{center}
\end{table}

We now turn to application of the results in \S\ref{s:eacrvav}.
It is easily shown that the conditions (\ref{eq:asbetatild})--(\ref{eq:asbetatild2}) hold for this model.
Using the same $N$ as above we obtain $\widetilde{\beta}_{0, N} \approx -0.417521221604055$, which coincides in all digits shown with $fl(\beta_0)$, evaluated as
$- \sum_{n=0}^N \bar{P}_n fl(\varphi_n)$.
Since $T_{p0} \approx 0.761$ and 
 $T_{10} \approx 1.117$, 
$T_{p0} - T_{10} \approx -0.356$, and hence, from 
Proposition \ref{pro:fbnbetasigma2phi}(a) one can argue heuristically that
\[
| E_{\textup{abs}}(\widetilde{\beta}_{0, N})| \approx |(T_{p0} - T_{10}) E_{\textup{abs}}(\widehat{\zeta})| \approx \widehat{\zeta} \, |(T_{p0} - T_{10}) E_{\textup{rel}}(\widehat{\zeta})| 
< 
0.143 \, u,
\]
and hence $| E_{\textup{rel}}(\widetilde{\beta}_{0, N})| \approx | E_{\textup{abs}}(\widetilde{\beta}_{0, N})/\widetilde{\beta}_{0, N}| < 0.34 \, u$. 

As for the asymptotic variance, we have
$\widetilde{\sigma}_N^2 \approx 0.589053281069282$, which matches $fl(\sigma^2)$, evaluated as 
$2 \sum_{n=0}^N \lambda_n p_n fl(\varphi_n)^2$, in all 15 digits.  
Furthermore, using 
Proposition \ref{pro:fbnbetasigma2phi}(b) and $\beta_1 \approx 0.025$, one can argue  that
\[
|E_{\textup{abs}}(\widetilde{\sigma}_N^2)| \approx  4 |\beta_{1} E_{\textup{abs}}(\widehat{\zeta})|  
 \approx 
4 \widehat{\zeta} \, |\beta_{1} E_{\textup{rel}}(\widehat{\zeta})| < 0.04 \, u,
\]
and hence $| E_{\textup{rel}}(\widetilde{\sigma}_N^2)| \approx | E_{\textup{abs}}(\widetilde{\sigma}_N^2)/\widetilde{\sigma}_N^2| < 0.07 \, u$. 

\section{Conclusions}
\label{s:cr}
While there is extensive work on the numerical instability analysis of linear recurrences, to date there is a dearth of research on its application to recurrences arising in applied probability, such as the Poisson equation considered herein.
In this paper, the rich structure of this equation has been exploited to develop an error analysis elucidating the instability phenomenon in its numerical solution, as well as a means of overcoming it.
It would be interesting to extend these results beyond the present scope to general continuous-time Markov chains. 

\section{Acknowledgments}
\label{s:ack}
The author thanks two anonymous reviewers for constructive comments that helped improved the paper.

\appendix
\renewcommand{\thesection}{\Alph{section}}

\section{Groundwork for the proof of Theorem \ref{the:pofzi}}
\label{a:pofzi}
This appendix lays the groundwork for the proof of Theorem \ref{the:pofzi} given in \S\ref{s:pofzi}.
Recall from \S \ref{s:intro} that we write
$d_n \triangleq \mu_n - \lambda_n$ and  $\Delta x_n \triangleq x_n - x_{n-1}$.

We use below the following  identities: for  $a, b, c, d \in \mathbb{R}$ with $b, d \neq 0$,  
\begin{equation}
\label{eq:ei1}
\frac{a+c}{b+d} = \frac{a}{b} + \frac{d}{b+d} \bigg(\frac{c}{d} - \frac{a}{b}\bigg)
\end{equation}
and
\begin{equation}
\label{eq:ei2}
\frac{c}{d} + \frac{b}{b-d} \bigg(\frac{a}{b} - \frac{c}{d}\bigg)= \frac{a-c}{b-d} = \frac{a}{b} + \frac{d}{b-d} \bigg(\frac{a}{b} - \frac{c}{d}\bigg).
\end{equation}
We will further use the following inequalities: if $b, d > 0$,
\begin{equation}
\label{eq:ein1}
\frac{a}{b} \leqslant \frac{a+c}{b+d} \leqslant \frac{c}{d}   \quad \textup{if and only if} \quad \frac{a}{b} \leqslant \frac{c}{d},   
\end{equation}
which is the classic \emph{mediant inequality}. Further, if $b > d > 0$, 
\begin{equation}
\label{eq:ein2}
\frac{a}{b} \leqslant \frac{a-c}{b-d}   \quad \textup{if and only if} \quad \frac{c}{d} \leqslant  \frac{a}{b}.
\end{equation}
Note that, in both (\ref{eq:ei2}) and (\ref{eq:ein2}), the leftmost inequalities are strict if and only if the rightmost inequalities are strict.

\begin{lemma}
\label{lma:ziinc} Under Assumption \textup{\ref{ass:hmulambda}(ii.a),} 
\begin{itemize}
\item[\textup{(a)}]
$Z_n$ is nondecreasing$;$
\item[\textup{(b)}] $\varphi_n$ is nonnegative$.$
\end{itemize}
\end{lemma}
\begin{proof}
(a)
Let $n \geqslant 1$. Then, 
using  (\ref{eq:Zirel}), (\ref{eq:Pi}), (\ref{eq:Hi}),  and (\ref{eq:ein1}), we obtain
\begin{equation}
\label{eq:zimzim1}
\Delta Z_n =  
\frac{\mu_n H_{n-1}^+ + c_n}{\mu_n T_{n-1}^+ + 1} - \frac{H_{n-1}^+}{T_{n-1}^+} = \frac{c_n - Z_{n-1}}{\mu_n T_{n-1}^+ + 1}.
\end{equation}
It now follows that $\Delta Z_n \geqslant 0$ since, by (\ref{eq:PC}), (\ref{eq:Zirel}), and Assumption \ref{ass:hmulambda}(ii.a),
\[
Z_{n-1} = \frac{\sum_{j=0}^{n-1} c_{j} p_{j}}{\sum_{j=0}^{n-1} p_{j}} \leqslant c_n.
\]

(b) The result follows by part (a), Corollary \ref{cor:phiexp}(a) and (\ref{eq:PCZ1is}).
\end{proof}

\begin{lemma}
\label{lma:lambDHPsti} \hspace{1in}
\begin{itemize}
\item[\textup{(a)}] 
\begin{align*}
\lambda_n \Delta H_n^+  & =  
\begin{cases}
\Delta c_1 + H_0^+ \Delta d_1,& \quad n = 1 \\
\Delta c_n + \mu_{n-1} \Delta H_{n-1}^+   + H_{n-1}^+ \Delta d_n, & \quad n \geqslant 2
\end{cases} \\
& = \Delta c_n + (\mu_{n-1} + \Delta d_n) \Delta H_{n-1}^+   + H_{n-2}^+ \Delta d_n, \qquad n \geqslant 2;
\end{align*}
\item[\textup{(b)}] 
\begin{align*}
\lambda_n \Delta T_n^+ & = 
\begin{cases}
T_{0} \Delta d_1, & \quad n = 1 \\
\mu_{n-1} \Delta T_{n-1}^+  + T_{n-1}^+ \Delta d_n, & \quad n \geqslant 2
\end{cases} \\
& = (\mu_{n-1} + \Delta d_n) \Delta T_{n-1}^+  + T_{n-2}^+ \Delta d_n, \quad n \geqslant 2.
\end{align*}
\end{itemize}
\end{lemma}
\begin{proof} 
(a) 
From (\ref{eq:Hi}) we obtain
\begin{equation}
\label{eq:hdHstpii}
\lambda_n \Delta H_n^+  = c_n + d_n H_{n-1}^+, \quad n \geqslant 1,
\end{equation}
For $n = 1$, using (\ref{eq:hdHstpii}) and $c_{0} + d_{0} H_0^+ = 0$ (cf.\ (\ref{eq:Hi})), we obtain
\begin{align*}
\lambda_1 \Delta H_1^+ & = c_1 + d_1 H_0^+ =  c_1 + H_0^+ \Delta d_1 + d_{0} H_0^+   = \Delta c_1 + H_0^+ \Delta d_1.
\end{align*}

For $n \geqslant 2$, using twice (\ref{eq:hdHstpii}) yields
\begin{align*}
\lambda_n \Delta H_n^+  & = c_n + d_n H_{n-1}^+ = 
\Delta c_n +  H_{n-1}^+ \Delta d_n + c_{n-1} +  d_{n-1} H_{n-1}^+
\\
& = \Delta c_n +  H_{n-1}^+ \Delta d_n + c_{n-1} +  d_{n-1} H_{n-2}^+ + d_{n-1} \Delta H_{n-1}^+ \\
& = \Delta c_n +  H_{n-1}^+ \Delta d_n + \lambda_{n-1} \Delta H_{n-1}^+ + d_{n-1} \Delta H_{n-1}^+ \\
& = \Delta c_n + \mu_{n-1} \Delta H_{n-1}^+ +  H_{n-1}^+ \Delta d_n.
\end{align*}

(b) From (\ref{eq:Pi}) we readily obtain
\begin{equation}
\label{eq:hdPstpii}
\lambda_n \Delta T_n^+  = 1 + d_n T_{n-1}^+, \quad n \geqslant 1.
\end{equation}
This part follows 
 as part (a) using (\ref{eq:hdPstpii}) and $1 + d_{0} T_{0} = 0$ (cf.\ (\ref{eq:Pi})).
\end{proof}

\begin{lemma}
\label{lma:PHstarinc} Under Assumption \textup{\ref{ass:hmulambda}(i.a, ii.a),} 
\begin{itemize}
\item[\textup{(a)}] 
$H_n^+$ is nondecreasing$;$
\item[\textup{(b)}]
$T_n^+$ is increasing$.$
\end{itemize}
\end{lemma}
\begin{proof}
(a) 
This part follows immediately 
by induction using Lemma \ref{lma:lambDHPsti}(a) and Assumption \ref{ass:hmulambda}(i.a, ii.a), which yields $H_{n-1}^+ \geqslant 0$ and $\Delta H_n^+ \geqslant 0$ for $n \geqslant 1$.

(b) This part follows similarly using Lemma \ref{lma:lambDHPsti}(b) and Assumption \ref{ass:hmulambda}(i.a),  which yields $\Delta T_n^+ > 0$ for $n \geqslant 1$.
\end{proof}

\begin{lemma}
\label{lma:difDeltaziinc} Under Assumption $\ref{ass:hmulambda}($\textup{i.a, ii.a}$),$
$Z_n \leqslant {\Delta H_n^+}/{\Delta T_n^+}.$ 
\end{lemma}
\begin{proof}
We have, using (\ref{eq:ei2}), (\ref{eq:ein2}), and Lemmas \ref{lma:ziinc} and \ref{lma:PHstarinc},
\begin{equation}
\label{eq:deltaHPsteq}
\frac{\Delta H_n^+}{\Delta T_n^+} = 
 \frac{H_n^+}{T_n^+} + \frac{T_{n-1}^+}{\Delta T_n^+} \bigg(\frac{H_n^+}{T_n^+} - \frac{H_{n-1}^+}{T_{n-1}^+}\bigg) \geqslant
 \frac{H_n^+}{T_n^+} = Z_n.
\end{equation}
\end{proof}

\begin{lemma}
\label{lma:DeltaHPHPineq} Under Assumption $\ref{ass:hmulambda}($\textup{i, ii)}$,$
\begin{equation}
\label{eq:DeltaHPHPineq}
 \bigg(\frac{\Delta H_n^+}{\Delta T_n^+}-\frac{H_{n-1}^+}{T_{n-1}^+}\bigg) \Delta d_n \leqslant 
\frac{\Delta c_n}{T_{n-1}^+}, \quad n \geqslant 1 \; \textup{(with equality for $n = 1.)$}
\end{equation}
\end{lemma}
\begin{proof} 
We prove the result by induction.
For $n = 1$, using Lemmas \ref{lma:lambDHPsti} and \ref{lma:PHstarinc}(b), and Assumption \ref{ass:hmulambda}(i.a), it follows that (\ref{eq:DeltaHPHPineq}) holds with equality:

\begin{equation}
\label{eq:dhs1dps1eq}
\begin{split}
\frac{\Delta H_1^+}{\Delta T_1^+} & = 
\frac{H_0^+ \Delta d_1 + \Delta c_1 }{T_{0} \Delta d_1} =
\frac{H_0^+}{T_{0}} + 
\frac{\Delta c_1}{T_{0} \Delta d_1}.
\end{split}
\end{equation}

Now, suppose that (\ref{eq:DeltaHPHPineq}) holds for some $n \geqslant 1$. 
If $\Delta d_{n+1} = 0$, then it trivially holds for $n+1$, since $\Delta c_{n+1} \geqslant 0$ by Assumption \ref{ass:hmulambda}(ii.a).

So consider the case $\Delta d_{n+1} \neq 0$. Then 
$0 < \Delta d_{n+1} \leqslant \Delta d_n$ by Assumption \ref{ass:hmulambda}(i).
Using Lemmas \ref{lma:lambDHPsti} and \ref{lma:PHstarinc}(b), and (\ref{eq:ei1}), we can write

\begin{align*}
\frac{\Delta H_{n+1}}{\Delta T_{n+1}} & = \frac{H_n^+ \Delta d_{n+1} + \mu_n \Delta H_n^+ + \Delta c_{n+1}}{T_n^+ \Delta d_{n+1} + \mu_n \Delta T_n^+} \\
& =
\frac{H_n^+}{T_n^+} + 
\frac{\mu_n \Delta T_n^+}{\lambda_{n+1} \Delta T_{n+1}}\bigg(\frac{\mu_n \Delta H_n^+ + \Delta c_{n+1}}{\mu_n \Delta T_n^+} - \frac{H_n^+}{T_n^+}\bigg).
\end{align*}

Hence, we can reformulate the required result that (\ref{eq:DeltaHPHPineq}) holds for $i + 1$ as
\begin{equation}
\label{eq:dhs1dps1eqip1}
\frac{\mu_n \Delta T_n^+}{\lambda_{n+1} \Delta T_{n+1}}\bigg(\frac{\mu_n \Delta H_n^+ + \Delta c_{n+1}}{\mu_n \Delta T_n^+} - \frac{H_n^+}{T_n^+}\bigg) \leqslant \frac{\Delta c_{n+1}}{T_n^+ \Delta d_{n+1}}.
\end{equation}
In turn, we can reformulate the latter inequality as follows:
\[
\frac{T_n^+ \Delta d_{n+1}}{\lambda_{n+1} \Delta T_{n+1}}\bigg(\mu_n \Delta H_n^+ + \Delta c_{n+1} - \mu_n \Delta T_n^+ \frac{H_n^+}{T_n^+}\bigg) \leqslant \Delta c_{n+1},
\]
i.e., 
\[
\frac{T_n^+ \Delta d_{n+1}}{\lambda_{n+1} \Delta T_{n+1}}\bigg(\Delta c_{n+1} + \mu_n \Delta T_n^+ \bigg(\frac{\Delta H_n^+}{\Delta T_n^+}  - \frac{H_n^+}{T_n^+}\bigg)\bigg) \leqslant \Delta c_{n+1},
\]
i.e., using again Lemma \ref{lma:lambDHPsti}(b),
\[
\frac{T_n^+ \Delta d_{n+1} \mu_n \Delta T_n^+}{\lambda_{n+1} \Delta T_{n+1}} \bigg(\frac{\Delta H_n^+}{\Delta T_n^+}  - \frac{H_n^+}{T_n^+}\bigg) \leqslant 
\frac{\mu_n \Delta T_n^+}{\lambda_{n+1} \Delta T_{n+1}}\Delta c_{n+1},
\]
i.e., 
\begin{equation}
\label{eq:dHiPiHiPi}
 \frac{\Delta H_n^+}{\Delta T_n^+}  - \frac{H_n^+}{T_n^+} \leqslant 
\frac{\Delta c_{n+1}}{T_n^+ \Delta d_{n+1}}.
\end{equation}

To prove (\ref{eq:dHiPiHiPi}), we write, using (\ref{eq:ei1}),
\begin{align*}
\frac{H_n^+}{T_n^+} + \frac{\Delta c_{n+1}}{T_n^+ \Delta d_{n+1}} & = 
\frac{H_n^+ \Delta d_{n+1} + \Delta c_{n+1}}{T_n^+ \Delta d_{n+1}} \\
& = 
\frac{(\Delta H_n^+) \Delta d_{n+1} + H_{n-1}^+ \Delta d_{n+1} + \Delta c_{n+1}}{(\Delta T_n^+) \Delta d_{n+1} + T_{n-1}^+ \Delta d_{n+1}} \\
& =\frac{\Delta H_n^+}{\Delta T_n^+} + \frac{T_{n-1}^+}{T_n^+}
\bigg(\frac{H_{n-1}^+}{T_{n-1}^+} + \frac{\Delta c_{n+1}}{T_{n-1}^+ \Delta d_{n+1}} - \frac{\Delta H_n^+}{\Delta T_n^+}\bigg),
\end{align*}
whence
\begin{align*}
\frac{H_n^+}{T_n^+} + \frac{\Delta c_{n+1}}{T_n^+ \Delta d_{n+1}} - \frac{\Delta H_n^+}{\Delta T_n^+} 
& = \frac{T_{n-1}^+}{T_n^+}
\bigg(\frac{H_{n-1}^+}{T_{n-1}^+} + \frac{\Delta c_{n+1}}{T_{n-1}^+ \Delta d_{n+1}} - \frac{\Delta H_n^+}{\Delta T_n^+}\bigg) \\
& \geqslant \frac{T_{n-1}^+}{T_n^+}
\bigg(\frac{H_{n-1}^+}{T_{n-1}^+} + \frac{\Delta c_n}{T_{n-1}^+ \Delta d_n} - \frac{\Delta H_n^+}{\Delta T_n^+}\bigg)  \geqslant 0,
\end{align*}
where the first and second inequalities follow by Assumption \ref{ass:hmulambda}(i.b, ii.b) and the induction hypothesis (\ref{eq:DeltaHPHPineq}), respectively. Therefore, (\ref{eq:dHiPiHiPi}) holds, and hence so does (\ref{eq:dhs1dps1eqip1}), which completes the induction proof.\end{proof}

\begin{lemma}
\label{lma:Deltaziinc} Under Assumption \textup{\ref{ass:hmulambda}(i, ii),}
$\displaystyle \frac{\Delta H_n^+}{\Delta T_n^+}$ is nondecreasing$.$ 
\end{lemma}
\begin{proof} 
Fix $n \geqslant 1$. We can write, using Lemma \ref{lma:lambDHPsti},
\begin{equation}
\label{eq:dhs2dps2eq}
\frac{\Delta H_{n+1}}{\Delta T_{n+1}} = 
\frac{(\mu_n + \Delta d_{n+1}) \Delta H_n^+ + \Delta c_{n+1} +
   H_{n-1}^+ \Delta d_{n+1}}{(\mu_n + \Delta d_{n+1}) \Delta T_n^+ +
   T_{n-1}^+ \Delta d_{n+1}}.
\end{equation}
We need to distinguish two cases. If  $\Delta d_{n+1} = 0$, the latter identity gives
\[
\frac{\Delta H_{n+1}}{\Delta T_{n+1}} = 
\frac{\mu_n  \Delta H_n^+ + \Delta c_{n+1}}{\mu_n \Delta T_n^+} = 
\frac{\Delta H_n^+ }{\Delta T_n^+} + 
\frac{\Delta c_{n+1}}{\mu_n \Delta T_n^+} \geqslant 
\frac{\Delta H_n^+ }{\Delta T_n^+},
\]
as required, using Lemma \ref{lma:PHstarinc}(b) and Assumption \ref{ass:hmulambda}(ii.a). 

If $\Delta d_{n+1} \neq 0$, it must be $\Delta d_{n+1} > 0$ by Assumption \ref{ass:hmulambda}(i.a). We have
\begin{align*} 
\frac{\Delta H_{n+1}}{\Delta T_{n+1}} & = 
\frac{\Delta H_n^+}{\Delta T_n^+} + 
\frac{T_{n-1}^+ \Delta d_{n+1}}{\lambda_{n+1} \Delta T_{n+1}} \bigg(\frac{\Delta c_{n+1}}{T_{n-1}^+ \Delta d_{n+1}}  - \bigg(\frac{\Delta H_n^+}{\Delta T_n^+}-\frac{H_{n-1}^+}{T_{n-1}^+}\bigg)\bigg) \\
& \geqslant 
\frac{\Delta H_n^+}{\Delta T_n^+} + 
\frac{T_{n-1}^+ \Delta d_{n+1}}{\lambda_{n+1} \Delta T_{n+1}}  \bigg(\frac{\Delta c_{n+1}}{T_{n-1}^+ \Delta d_{n+1}} - \frac{\Delta c_n}{T_{n-1}^+ \Delta d_n}\bigg)  \geqslant 
\frac{\Delta H_n^+}{\Delta T_n^+},
\end{align*}
using (\ref{eq:dhs2dps2eq}), (\ref{eq:ei1}), Lemmas \ref{lma:lambDHPsti}(b), \ref{lma:PHstarinc}(b) and \ref{lma:DeltaHPHPineq},    
and Assumption \ref{ass:hmulambda}(i.b, ii.b). 
\end{proof}

\begin{lemma}
\label{lma:limDeltaziinc} 
$\displaystyle \frac{\Delta H_n^+}{\Delta T_n^+} \to \zeta$ as $n \to \infty.$
\end{lemma}
\begin{proof}
We can write  
\begin{align*}
\frac{\Delta H_n^+}{\Delta T_n^+} - 
 \frac{H_n^+}{T_n^+} & = \frac{T_{n-1}^+}{\Delta T_n^+} \Delta Z_n   =  \frac{T_{n-1}^+}{\Delta T_n^+}
 \frac{1}{1 + \mu_n T_{n-1}^+} (c_n - Z_{n-1}) \\
 & = \frac{T_{n-1}^+}{\lambda_n \Delta T_n^+} \frac{c_n - Z_{n-1}}{T_n^+} = \frac{T_{n-1}^+}{1 + d_n T_{n-1}^+} \frac{c_n - Z_{n-1}}{T_n^+},
\end{align*}
where we have used in turn (\ref{eq:Zirel}), (\ref{eq:deltaHPsteq}), (\ref{eq:zimzim1}), (\ref{eq:Pi}), and (\ref{eq:hdPstpii}).

Now, we have
\begin{equation}
\label{eq:ktbounded}
\frac{T_{n-1}^+}{1 + d_n T_{n-1}^+} = 
\frac{1}{{1}/{T_{n-1}^+} + d_n} \to 
\frac{1}{{1}/{T_{\infty}^+} + d_{\infty}} < \infty \textup{  as  } n \to \infty,
\end{equation}
where $0< T_{\infty}^+ \leqslant \infty$ and
$0< d_{\infty} \leqslant \infty$ are the limits, possibly infinite, of $T_{n}^+$ and 
$d_n$ as $n \to \infty$. See Lemma \ref{lma:PHstarinc}(b) and Assumption \ref{ass:hmulambda}(i). Note that Assumption \ref{ass:hmulambda}(i.a) and ergodicity ensure that 
$d_{\infty} > 0.$ Otherwise, it would be $\mu_n \leqslant \lambda_n$ for all $n$, and the chain would not be ergodic.

Furthermore, using (\ref{eq:zimzim1}) and (\ref{eq:PCZ1is}) gives

\[
\frac{c_n - Z_{n-1}}{T_n^+} = \Delta Z_n \to 0 \textup{  as  } n \to \infty.
\]

Therefore we obtain, as required,
\[
\lim_{n \to \infty} \, \frac{\Delta H_n^+}{\Delta T_n^+} =  
 \lim_{n \to \infty} \, \frac{H_n^+}{T_n^+} = \zeta.
\]
\end{proof}



\begin{thebibliography}{10}
\expandafter\ifx\csname url\endcsname\relax
  \def\url#1{\texttt{#1}}\fi
\expandafter\ifx\csname urlprefix\endcsname\relax\def\urlprefix{URL }\fi
\expandafter\ifx\csname href\endcsname\relax
  \def\href#1#2{#2} \def\path#1{#1}\fi

\bibitem{ross19}
S.~M. Ross, Introduction to Probability Models, 12th Edition, Academic Press,
  New York, 2019.

\bibitem{stidh09}
S.~Stidham, Jr., Optimal Design of Queueing Systems, Chapman and Hall/CRC, Boca
  Raton, FL, 2009.

\bibitem{neveu71}
J.~Neveu, Recurrent {M}arkovian potential and {P}oisson equation, C. R. Acad.
  Sci., A 272 (1971) 1202--1204.

\bibitem{nummelin91}
E.~Nummelin, On the {P}oisson equation in the potential theory of a single
  kernel, Math. Scand. 68 (1991) 59--82.

\bibitem{glynnMeyn96}
P.~W. Glynn, S.~P. Meyn, A {L}iapounov bound for solutions of the {P}oisson
  equation, Ann. Probab. 24 (1996) 916--931.

\bibitem{meynTweed09}
S.~P. Meyn, R.~L. Tweedie, Markov Chains and Stochastic Stability, 2nd Edition,
  Cambridge University Press, Cambridge, UK, 2009.

\bibitem{guoHL09}
X.~Guo, O.~Hern\'andez-Lerma, Continuous-Time Markov Decision Processes: Theory
  and Applications, Springer, Heidelberg, Germany, 2009.

\bibitem{howard60}
R.~A. Howard, Dynamic {P}rogramming and {M}arkov {P}rocesses, Wiley, New York,
  1960.

\bibitem{blokSpieks17}
H.~Blok, F.~M. Spieksma, Structures of optimal policies in {MDP}s with
  unbounded jumps: The state of our art, in: R.~J. Boucherie, N.~M. van Dijk
  (Eds.), Markov Decision Processes in Practice, Springer, Cham, Switzerland,
  2017, Ch.~5, pp. 131--186.

\bibitem{meynTweed93}
S.~P. Meyn, R.~L. Tweedie, Stability of {M}arkovian processes {III}:
  {F}oster--{L}yapunov criteria for continuous-time processes, Adv. Appl.
  Probab. 25 (1993) 518--548.

\bibitem{nmcor19}
J.~Ni\~{n}o{-}Mora, Resource allocation and routing in parallel multi-server
  queues with abandonments for cloud profit maximization, Comput. Oper. Res.
  103 (2019) 221--236.

\bibitem{higham02}
N.~J. Higham, Accuracy and Stability of Numerical Algorithms, 2nd Edition,
  SIAM, Philadelphia, PA, 2002.

\bibitem{steckHend07}
S.~G. Steckley, S.~G. Henderson, The error in steady-state approximations for
  the time-dependent waiting time distribution, Stoch. Models 23 (2007)
  307--332.

\bibitem{junejaShah06}
S.~Juneja, P.~Shahabuddin, Rare-event simulation techniques: An introduction
  and recent advances, in: S.~G. Henderson, B.~L. Nelson (Eds.), Simulation,
  North Holland, Amsterdam, The Netherlands, 2006, Ch.~11, pp. 291--350.

\bibitem{krish90}
K.~R. Krishnan, Joining the right queue: A state-dependent decision rule, IEEE
  Trans. Automat. Control 35 (1990) 104--108.

\bibitem{bhulai06}
S.~Bhulai, On the value function of the {M/Cox$(r)$/$1$} queue, J. Appl.
  Probab. 43 (2006) 363--376.

\bibitem{hyytiaRight16}
E.~Hyyti\"a, R.~Righter, Routing jobs with deadlines to heterogeneous parallel
  servers, Oper. Res. Lett. 44 (2016) 507--513.

\bibitem{nmcor12}
J.~Ni\~{n}o{-}Mora, Admission and routing of soft real-time jobs to
  multiclusters: {D}esign and comparison of index policies, Comput. Oper. Res.
  39 (2012) 3431--3444.

\bibitem{nmejor12}
J.~Ni\~{n}o{-}Mora, Towards minimum loss job routing to parallel heterogeneous
  multiserver queues via index policies, European J. Oper. Res. 220 (2012)
  705--715.

\bibitem{tijms17}
H.~C. Tijms, One-step improvement ideas and computational aspects, in: R.~J.
  Boucherie, N.~M. van Dijk (Eds.), Markov Decision Processes in Practice,
  Springer, Cham, Switzerland, 2017, Ch.~1, pp. 3--32.

\bibitem{bhulai17}
S.~Bhulai, Value function approximation in complex queueing systems, in: R.~J.
  Boucherie, N.~M. van Dijk (Eds.), Markov Decision Processes in Practice,
  Springer, Cham, Switzerland, 2017, Ch.~2, pp. 33--62.

\bibitem{whittAFMP92}
W.~Whitt, Asymptotic formulas for {M}arkov processes with applications to
  simulation, Oper. Res. 40 (1992) 279--291.

\bibitem{niHend15}
E.~C. Ni, S.~G. Henderson, How hard are steady-state queueing simulations?, ACM
  Trans. Model. Comput. Simul. 25 (article 27) (2015) 1--21.

\bibitem{kyriDim16}
E.~G. Kyriakidis, T.~D. Dimitrakos, A semi-{M}arkov decision model for the
  optimal control of a simple immigration-birth--death process through the
  introduction of a predator, Comm. Stat. Theor. Meth. 45 (2016) 3860--3873.

\bibitem{nmmp02}
J.~Ni\~{n}o{-}Mora, Dynamic allocation indices for restless projects and
  queueing admission control: {A} polyhedral approach, Math. Program. 93 (2002)
  361--413.

\bibitem{gautschi61}
W.~Gautschi, Recursive computation of certain integrals, J. ACM 8 (1961)
  21--40.

\bibitem{nmnetgcoop14}
J.~Ni\~{n}o{-}Mora, Overcoming numerical instability in one-step policy
  improvement for admission and routing to queues with firm deadlines, in:
  Proceedings of the 7th International Conference on Network Games, Control and Optimization (NetGCooP 2014), IEEE, New
  York, 2017, pp. 127--134.

\bibitem{makoShw02}
A.~M. Makowski, A.~Shwartz, The {P}oisson equation for countable {M}arkov
  chains: {P}robabilistic methods and interpretations, in: E.~A. Feinberg,
  A.~Shwartz (Eds.), Handbook of Markov Decision Processes: Methods and
  Applications, Springer, New York, 2002, pp. 269--303.

\bibitem{bhulaiSpieks03}
S.~Bhulai, F.~M. Spieksma, On the uniqueness of solutions to the {P}oisson
  equations for average cost {M}arkov chains with unbounded cost functions,
  Math. Meth. Oper. Res. 58 (2003) 221--236.

\bibitem{dendetal13}
S.~Dendievel, G.~Latouche, Y.~Liu, Poisson’s equation for discrete-time
  quasi-birth-and-death processes, Performance Eval. 70 (2013) 564--577.

\bibitem{binietal16}
D.~Bini, S.~Dendievel, G.~Latouche, B.~Meini, General solution of the {P}oisson
  equation for quasi-birth-and-death processes, SIAM J. Appl. Math. 76 (2016)
  2397--2417.

\bibitem{gautschi67}
W.~Gautschi, Computational aspects of three-term recurrence relations, SIAM
  Rev. 9 (1967) 24--82.

\bibitem{olver67}
F.~W.~J. Olver, Numerical solution of second-order linear difference equations,
  J. Res. Natl. Bur. Stand. B Math. Sci. 71B (1967) 111--129.

\bibitem{amosBurg73}
D.~E. Amos, J.~W. Burgmeier, Computation with three-term, linear,
  nonhomogeneous recursion relations, SIAM Rev. 15 (1973) 335--351.

\bibitem{barrioetal03}
R.~Barrio, B.~Melendo, S.~Serrano, On the numerical evaluation of linear
  recurrences, J. Comput. Appl. Math. 150 (2003) 71--86.

\bibitem{karlMcG57}
S.~Karlin, J.~McGregor, The classification of birth and death processes, Trans.
  Am. Math. Soc. 86 (1957) 366--401.

\bibitem{keilson65}
J.~Keilson, A review of transient behavior in regular diffusion and birth--death
  processes. {P}art {II}, J. Appl. Probab. 2 (1965) 405--428.

\bibitem{feller59}
W.~Feller, The birth-and-death processes as diffusion processes, J. Math. Pures
  Appl. 38 (1959) 301--345.

\bibitem{CallKeil73II}
H.~Callaert, J.~Keilson, On exponential ergodicity and spectral structure for
  birth--death processes {II}, Stoch. Process. Appl. 1 (1973) 217--235.

\bibitem{CallKeil73I}
H.~Callaert, J.~Keilson, On exponential ergodicity and spectral structure for
  birth--death processes {I}, Stoch. Process. Appl. 1 (1973) 187--216.

\bibitem{coolenVanDoorn02}
P.~Coolen-Schrijner, E.~A. {van Doorn}, The deviation matrix of a
  continuous-time {M}arkov chain, Probab. Eng. Inform. Sci. 16 (2002) 351--366.

\bibitem{muresan08}
M.~Mure\c{s}an, A Concrete Approach to Classical Analysis, Springer, New York,
  2008.

\bibitem{meyn07}
S.~Meyn, Control Techniques for Complex Networks, Cambridge University Press,
  Cambridge, UK, 2007.

\bibitem{anckGaf62}
C.~J. Ancker, {Jr.}, A.~V. Gafarian, Queueing with impatient customers who
  leave at random, J. Ind. Eng. 13 (1962) 84--90.

\bibitem{NIST:DLMF}
{\it NIST Digital Library of Mathematical Functions}, http://dlmf.nist.gov/,
  Release 1.0.26 of 2020-03-15, {F}.~W.~J. Olver et al., eds.

\end{thebibliography}

\end{document}